\documentclass[11pt]{article}

\usepackage[latin1]{inputenc}\usepackage{epsfig}\usepackage{amsfonts}

\title{
  {\huge Fourier Theory on the Complex Plane II} \\
  Weak Convergence, Classification and \\
  Factorization of Singularities }

\author{
  \Large Jorge L. deLyra \\
  Department of Mathematical Physics \\
  Physics Institute \\
  University of São Paulo }

\date{April 27, 2015}

\newcommand{\ii}{\mbox{\boldmath$\imath$}}

\newcommand{\at}[2]{\left.\rule{0em}{3ex}\right[_{\,#1}^{\,#2}}

\newcommand{\e}[1]{\,{\rm e}^{#1}}

\newcommand{\ldot}{\mbox{\Large$\cdot$}\!}

\begin{document}\maketitle

\begin{abstract}
  \noindent
  The convergence of DP Fourier series which are neither strongly
  convergent nor strongly divergent is discussed in terms of the Taylor
  series of the corresponding inner analytic functions. These are the
  cases in which the maximum disk of convergence of the Taylor series of
  the inner analytic function is the open unit disk. An essentially
  complete classification, in terms of the singularity structure of the
  corresponding inner analytic functions, of the modes of convergence of a
  large class of DP Fourier series, is established. Given a weakly
  convergent Fourier series of a DP real function, it is shown how to
  generate from it other expressions involving trigonometric series, that
  converge to that same function, but with much better convergence
  characteristics. This is done by a procedure of factoring out the
  singularities of the corresponding inner analytic function, and works
  even for divergent Fourier series. This can be interpreted as a
  resummation technique, which is firmly anchored by the underlying
  analytic structure.
\end{abstract}

\section{Introduction}

In this paper we will discuss the question of the convergence of DP
Fourier series, in the light of the correspondence between FC pairs of DP
Fourier series and Taylor series of inner analytic functions, which was
established in a previous paper~\cite{FTotCPI}. This discussion was
started in that paper, and will be continued here in order to include the
more complex and difficult cases. We will assume that the reader is aware
of the contents of that paper, and will use without too much explanation
the concepts, the definitions and the notations established there. We will
limit ourselves here to the explanation of a few abbreviations, such as
DP, which stands for Definite Parity, and FC, which stands for Fourier
Conjugate, and to the restatement of the more basic concepts.

We refer as the {\em basic convergence theorem} of complex analysis to the
result that, if a complex power series around a point $z_{0}$ converges at
a point $z_{1}\neq z_{0}$, then it is convergent and absolutely convergent
in the interior of a disk centered at $z=z_{0}$ with its boundary passing
through $z_{1}$. In addition to this, it converges uniformly on any closed
set contained within this open disk. The concept of {\em inner analytic
  function} makes reference to a complex function that is analytic at
least on the open unit disk, assumes the value zero at $z=0$, and is the
analytic continuation of a real function on the interval $(-1,1)$ of the
real axis. The restriction of this complex function to the unit circle
results in a FC pair of DP real functions. The corresponding Taylor series
around $z=0$ converges at least on the open unit disk, assumes the value
zero at $z=0$, and has real coefficients. The restriction of this complex
power series to the unit circle results in a FC pair of DP Fourier series.

The concept of {\em Definite Parity} or DP Fourier series refers to a
Fourier series which has only the sine terms, or only the cosine terms,
without the constant term, and therefore has a definite parity on the
periodic interval $[-\pi,\pi]$. In the same way, DP real functions are
those that have definite parity on the periodic interval. The concept of
{\em Fourier Conjugate} or FC trigonometric series refers to a conjugate
series that is built from a given DP trigonometric series by the exchange
of cosines by sines, or vice-versa. The concept of FC real functions
refers to the real functions that FC Fourier series converge to, or more
generally to the pair of real functions that generates the common Fourier
coefficients of a FC pair of DP trigonometric series, even if the series
do not converge.

The concepts of the {\em degrees of hardness} and of the {\em degrees of
  softness} of singularities in the complex plane, which were introduced
in the previous paper~\cite{FTotCPI}, will be discussed and defined in
more precise terms in this paper. The basic concept is that a soft
singularity is one where the complex function is still well defined,
although not analytic, when one takes the limit to the singular point. On
the other hand, a hard singularity is one where the complex function
diverges to infinity or does not exist when one takes that limit. A
borderline hard singularity is the least hard type of hard singularity,
while a borderline soft singularity is the least soft type of soft
singularity.

We will consider, then, the convergence of DP Fourier series and of the
corresponding complex power series. For organizational reasons this
discussion must be separated into several parts. In the previous
paper~\cite{FTotCPI} we tackled the extreme cases which we qualify as
those of {\em strong divergence} and {\em strong convergence}, for which
results can be established in full generality. In this paper we will deal
with the remaining cases, which we will qualify as those of {\em weak
  convergence} and {\em very weak convergence}. In the case of very weak
convergence we will be able to present definite results only for a certain
class of DP Fourier series.

\section{Weak Convergence}

In the previous paper~\cite{FTotCPI} we discussed the issue of the
convergence of FC pairs of DP Fourier series such as

\noindent
\begin{eqnarray*}
  S_{\rm c}
  & = &
  \sum_{k=1}^{\infty}
  a_{k}\cos(k\theta),
  \\
  S_{\rm s}
  & = &
  \sum_{k=1}^{\infty}
  a_{k}\sin(k\theta),
\end{eqnarray*}

\noindent
with real Fourier coefficients $a_{k}$, for $\theta$ in the periodic
interval $[-\pi,\pi]$, and of the corresponding complex power series

\noindent
\begin{eqnarray*}
  S_{v}
  & = &
  \sum_{k=1}^{\infty}
  a_{k}v^{k},
  \\
  S_{z}
  & = &
  \sum_{k=1}^{\infty}
  a_{k}z^{k},
\end{eqnarray*}

\noindent
where $v=\exp(\ii\theta)$ and $z=\rho v$, in two extreme cases, those of
strong divergence and strong convergence. We established that, in these
extreme cases, either the convergence/divergence of $S_{z}$ or the
existence/absence of singularities of the analytic function $w(z)$ that
$S_{z}$ converges to can be used to determine the convergence or
divergence of $S_{v}$ and of the corresponding DP Fourier series over the
whole unit circle of the complex plane.

Specifically, we established that the divergence of $S_{z}$, or the
existence of a singularity of $w(z)$, at any single point strictly within
the open unit disk, implies that $S_{v}$ diverges everywhere on the unit
circle, and that the corresponding FC pair of DP Fourier series diverge
almost everywhere on the periodic interval. This is what we call strong
divergence. To complement this we established that the convergence of
$S_{z}$ at a single point strictly outside the closed unit disk, or the
absence of singularities of $w(z)$ on that disk, implies that $S_{v}$
converges to a $C^{\infty}$ function everywhere on the unit circle, and
that the corresponding FC pair of DP Fourier series converge to
$C^{\infty}$ real functions everywhere on the periodic interval. This is
what we call strong convergence.

If the series $S_{z}$ constructed from an FC pair of DP trigonometric
series converges on at least one point on the unit circle, but does not
converge at any points strictly outside the closed unit disk, then the
situation is significantly more complex. We will see, however, that in
essence it can still be completely determined, with some limitations on
the set of series involved when establishing the results. We will refer to
this situation as the case of {\em weak convergence}. In this situation
the maximum disk of convergence of the complex power series $S_{z}$ is the
open unit disk, and therefore this series converges strongly to an
analytic function $w(z)$ strictly within that disk. As is well known, it
can be shown that in this case the function $w(z)$ must have at least one
singularity on the unit circle.

In terms of the convergence of the series, however, the situation on the
unit circle remains undefined. The series $S_{z}$ may converge or diverge
at various points on that circle. The function $w(z)$ cannot be analytic
on the whole unit circle, since there must be at least one point on that
circle where it has a singularity, in the sense that it is not analytic at
that point. This singularity can be of various types, and does not
necessarily mean that the function $w(z)$ is not defined at its location.
In fact, it is still possible for $w(z)$ to exist everywhere and to be
analytic almost everywhere on the unit circle. No general convergence
theorem is available for this case, and the analysis must be based on the
behavior of the coefficients $a_{k}$ of the series. It is, therefore,
significantly harder to obtain definite results.

While the case of strong divergence is characterized by the fact that the
coefficients $a_{k}$ typically diverge exponentially with $k$, and the
case of strong convergence by the fact that they typically go to zero
exponentially with $k$, in the case of weak convergence the typical
behavior is that the coefficients go to zero as a negative power of $k$.
In most of this paper we will in fact limit the discussion to series in
which the coefficients go to zero as a negative but not necessarily
integer power of $k$.

\subsection{Absolute and Uniform Convergence}

The simplest case of what we call here weak convergence is that in which
there is absolute and uniform convergence at the unit circle to a function
which, although necessarily continuous, is not $C^{\infty}$. If the
coefficients $a_{k}$ are such that the series $S_{z}$ is absolutely
convergent on a single point $z_{1}$ of the unit circle, that is, if the
series of absolute values

\begin{displaymath}
  \sum_{k=1}^{\infty}
  |a_{k}||z_{1}|^{k}
  =
  \sum_{k=1}^{\infty}
  |a_{k}|,
\end{displaymath}

\noindent
where $|z_{1}|=1$, is convergent, then the series $S_{v}$ is absolutely
convergent on the whole unit circle, given that the criterion of
convergence is clearly independent of $\theta$ in this case. For the same
reason, the series in uniformly convergent and thus converges to a
continuous function over the whole unit circle. Therefore the
corresponding FC Fourier series are also absolutely and uniformly
convergent to continuous functions $f_{\rm c}(\theta)$ and $f_{\rm
  s}(\theta)$ over their whole domain.

In this case, although the function $w(z)$ must have at least one
singularity on the unit circle, since the series $S(z)$ converges
everywhere that function is still well defined over the whole unit circle,
and by Abel's theorem~\cite{AbelTheo} the $\rho\to 1$ limit of the
function $w(z)$ from the interior of the unit disk to the unit circle
exists at all points of that circle. It follows that the singularity or
singularities of $w(z)$ at the unit circle must not involve divergences to
infinity. We will classify such singularities, at which the complex
function is still well-defined, although non-analytic, as {\em soft
  singularities}. Singularities at which the complex function diverges to
infinity or cannot be defined at all will be called {\em hard
  singularities}.

Since they are convergent everywhere, in this case all DP trigonometric
series are DP Fourier series of the continuous function obtained on the
unit circle by the $\rho\to 1$ limit of the real and imaginary parts of
$w(z)$ from within the unit disk. In short, we have established that
weakly-convergent FC pairs of DP Fourier series that converge absolutely
and thus uniformly must correspond to inner analytic functions that have
only soft singularities on the unit circle. In fact, this would be true of
any FC pair of DP Fourier series that simply converges weakly everywhere
over the unit circle.

Most of the time this situation can be determined fairly easily in terms
of the coefficients $a_{k}$ of the series. If there is a positive real
constant $A$, a strictly positive number $\varepsilon$ and an integer
$k_{m}$ such that for $k>k_{m}$ the absolute values of the coefficients
can be bounded as

\begin{displaymath}
  |a_{k}|
  \leq
  \frac{A}{k^{1+\varepsilon}},
\end{displaymath}

\noindent
then the asymptotic part of the series of the absolute values of the terms
of $S_{v}$, which is a sum of positive terms and thus increases
monotonically, can be bounded from above by a convergent asymptotic
integral, as shown in Section~\ref{APPevalconv} of
Appendix~\ref{APPproofs}, and thus that series converges. It follows that
the $S_{v}$ series is absolutely and uniformly convergent on the whole
unit circle, and therefore so are the two corresponding DP Fourier series.

If the coefficients tend to zero as a power when $k\to\infty$, then it is
easy to see that the limiting function defined by the series cannot be
$C^{\infty}$. Since every term-wise derivative of the Fourier series with
respect to $\theta$ adds a factor of $k$ to the coefficients, a decay to
zero as $1/k^{n+1+\varepsilon}$ with $0<\varepsilon\leq 1$ guarantees that
we may take only up to $n$ term-by-term derivatives and still end up with
an absolutely and uniformly convergent series. With $n+1$ differentiations
the resulting series might still converge almost everywhere, but that is
not certain, and typically at least one of the two DP series in the FC
pair will diverge somewhere. With $n+2$ differentiations the series
$S_{v}$ is sure to diverge everywhere, and the two DP series in the FC
pair are sure to diverge almost everywhere, since the coefficients no
longer go to zero as $k\to\infty$.

It is easy to see that the points where one of the limiting functions
$f_{\rm c}(\theta)$ and $f_{\rm s}(\theta)$ is not differentiable, or
those at which there is a singularity in one of its derivatives, must be
points where $w(z)$ is singular. At any point of the unit circle where
$w(z)$ is analytic it not only is continuous but also infinitely
differentiable. We see therefore that the set of points of the unit circle
where $f_{\rm c}(\theta)$ and $f_{\rm s}(\theta)$ or their derivatives of
any order have singularities of any kind must be exactly the same set of
points where $w(z)$ has singularities on that circle. In our present case,
since the functions must be continuous, these singularities can only be
points of non-differentiability of the functions or points where some of
their higher-order derivatives do not exist.

We conclude therefore that the cases in which the convergence is not
strong but is still absolute and uniform are easily characterized. The
really difficult cases regarding convergence are, therefore, those in
which the series are not absolutely or uniformly convergent on the unit
circle. These are cases in which $|a_{k}|$ goes to zero as $1/k$ or slower
as $k\to\infty$, as shown in Section~\ref{APPevalconv} of
Appendix~\ref{APPproofs}. We will refer to these cases as those of {\em
  very weak convergence}. In this case the Fourier series can converge to
discontinuous functions, and $S_{z}$ will typically diverge at some points
of the unit circle, at which $w(z)$ has divergent limits and therefore
hard singularities. In order to discuss this case we must first establish
a few simple preliminary facts, leading to a classification of
singularities according to their severity.

\section{Classification of Singularities}

Let us establish a general classification of the singularities of inner
analytic functions on the unit circle. While it is possible that the
scheme of classification that we describe here may have its uses in more
general settings, for definiteness we consider here only the case of inner
analytic functions that have one or more singularities on the unit circle.
We also limit our attention to only those singular points, and ignore any
singularities that may exist strictly outside the closed unit disk. In
order to do this we must first establish some preliminary facts.

To start with, let us show that the operation of logarithmic
differentiation stays within the set of inner analytic functions. Let us
recall that, if $w(z)$ is an inner analytic function, then it has the
properties that it is analytic on the open unit disk, that it is the
analytic continuation of a real function on the real interval $(-1,1)$,
and that $w(0)=0$. These last two are consequences of the fact that its
Taylor series has the form

\begin{displaymath}
  w(z)
  =
  \sum_{k=1}^{\infty}
  a_{k}z^{k},
\end{displaymath}

\noindent
with real $a_{k}$. We define the logarithmic derivative of $w(z)$ as

\begin{displaymath}
  w^{\ldot}(z)
  =
  z\,
  \frac{dw(z)}{dz},
\end{displaymath}

\noindent
which also establishes our notation for it. We might also write the first
logarithmic derivative as $w^{1\!\ldot}(z)$. Let us show that the
logarithmic derivative of an inner analytic function is another inner
analytic function. First, since the derivative of an analytic function is
analytic in exactly the same domain as that function, and since the
identity function $z$ is analytic in the whole complex plane, it follows
that if $w(z)$ is analytic on the open unit disk, then so is
$w^{\ldot}(z)$. Second, if we calculate $w^{\ldot}(z)$ using the series
representation of $w(z)$, which converges in the open unit disk, since
$w(z)$ is analytic there, we get

\begin{displaymath}
  w^{\ldot}(z)
  =
  z
  \sum_{k=1}^{\infty}
  ka_{k}z^{k-1},
\end{displaymath}

\noindent
since a convergent power series can always be differentiated term-by-term.
It follows that, if the coefficients $a_{k}$ are real, then so are the new
coefficients $ka_{k}$, so that the coefficients of the Taylor series of
$w^{\ldot}(z)$ are real, and hence it too is the analytic continuation of
a real function on the real interval $(-1,1)$. Lastly, due to the extra
factor of $z$ we have that $w^{\ldot}(0)=0$. Hence, the logarithmic
derivative $w^{\ldot}(z)$ is an inner analytic function.

Next, let us define the concept of logarithmic integration. This is the
inverse operation to logarithmic differentiation, with the understanding
that we always choose the value zero for the integration constant. The
logarithmic primitive of $w(z)$ may be defined as

\begin{displaymath}
  w^{-1\!\ldot}(z)
  =
  \int_{0}^{z}dz'\,
  \frac{1}{z'}\,
  w(z').
\end{displaymath}

\noindent
Let us now show that the operation of logarithmic integration also stays
within the set of inner analytic functions. Note that since $w(0)=0$ the
integrand is in fact analytic on the open unit disk, if we define it at
$z=0$ by continuity, and therefore the path of integration from $0$ to $z$
is irrelevant, so long as it is contained within that disk. It is easier
to see this using the series representation of $w(z)$, which converges in
the open unit disk, since $w(z)$ is analytic there,

\begin{displaymath}
  w^{-1\!\ldot}(z)
  =
  \int_{0}^{z}dz'\,
  \sum_{k=1}^{\infty}
  a_{k}z^{\prime(k-1)}.
\end{displaymath}

\noindent
It is clear now that the integrand is a power series which converges
within the open unit disk, and thus converges to an analytic function
there. Since the primitive of an analytic function is analytic in exactly
the same domain as that function, it follows that the logarithmic
primitive $w^{-1\!\ldot}(z)$ is analytic on the open unit disk. If we now
execute the integration using the series representation, we get

\begin{displaymath}
  w^{-1\!\ldot}(z)
  =
  \sum_{k=1}^{\infty}
  \frac{a_{k}}{k}\,
  z^{k},
\end{displaymath}

\noindent
since a convergent power series can always be integrated term-by-term.
Since the coefficients $a_{k}$ are real, so are the new coefficients
$a_{k}/k$, and therefore $w^{-1\!\ldot}(z)$ is the analytic continuation
of a real function on the real interval $(-1,1)$. Besides, one can see
explicitly that $w^{-1\!\ldot}(0)=0$. Therefore, $w^{-1\!\ldot}(z)$ is an
inner analytic function. It is now easy to see also that the logarithmic
derivative of this primitive gives us back $w(z)$.

Finally, let us show that the derivative of $w(z)$ with respect to
$\theta$ is given by the logarithmic derivative of $w(z)$. Since we have
that $z=\rho\exp(\ii\theta)$, we may at once write that

\noindent
\begin{eqnarray*}
  \frac{dw(z)}{d\theta}
  & = &
  \frac{dz}{d\theta}\,
  \frac{dw(z)}{dz}
  \\
  & = &
  \ii z\,
  \frac{dw(z)}{dz}
  \\
  & = &
  \ii w^{\ldot}(z).
\end{eqnarray*}

\noindent
In the limit $\rho\to 1$, where and when that limit exists, the derivative
of $w(z)$ with respect to $\theta$ becomes the derivatives of the limiting
functions $f_{\rm c}(\theta)$ and $f_{\rm s}(\theta)$ of the FC pair of DP
Fourier series. The factor of $\ii$ effects the interchange of real and
imaginary parts, and the change of sign, that are consequences of the
differentiation of the trigonometric functions.

We see therefore that, given an inner-analytic function and its set of
singularities on the unit circle, as well as the corresponding FC pair of
DP Fourier series, we may at once define a whole infinite chain of
inner-analytic functions and corresponding DP Fourier series, running by
differentiation to one side and by integration to the other, indefinitely
in both directions. We will name this an {\em integral-differential
  chain}. We are now ready to give the complete formal definition of the
proposed classification of the singularities of inner analytic functions
$w(z)$ on the unit circle. Let $z_{1}$ be a point on the unit circle. We
start with the very basic classification which was already mentioned.

\begin{itemize}

\item A singularity of $w(z)$ at $z_{1}$ is a {\em soft singularity} if
  the limit of $w(z)$ from within the unit disk to that point exists and
  is finite.

\item A singularity of $w(z)$ at $z_{1}$ is a {\em hard singularity} if
  the limit of $w(z)$ from within the unit disk to that point does not
  exist, or is infinite.

\end{itemize}

\noindent
Next we establish a gradation of the concepts of hardness and softness of
the singularities of $w(z)$. To each singular point $z_{1}$ we attach an
integer giving either its {\em degree of hardness} or its {\em degree of
  softness}. In order to do this the following definitions are adopted.

\begin{itemize}

\item A single logarithmic integration of $w(z)$, that does not change the
  hard/soft character of a singularity, increases the degree of softness
  of that singularity by one, if it is soft, or decreases the degree of
  hardness of that singularity by one, if it is hard.

\item A single logarithmic differentiation of $w(z)$, that does not change
  the hard/soft character of a singularity, increases the degree of
  hardness of that singularity by one, if it is hard, or decreases the
  degree of softness of that singularity by one, if it is soft.

\item Given a soft singularity at $z_{1}$, if $w(z)$ can be
  logarithmically differentiated indefinitely without that singularity
  ever becoming hard, then we say that it is an {\em infinitely soft
    singularity}. This means that its degree of softness is infinite.

\item Given a hard singularity at $z_{1}$, if $w(z)$ can be
  logarithmically integrated indefinitely without that singularity ever
  becoming soft, then we say that it is an {\em infinitely hard
    singularity}. This means that its degree of hardness is infinite.

\item A singularity of $w(z)$ at $z_{1}$ is a {\em borderline soft} one if
  it is a soft singularity and a single logarithmic differentiation of
  $w(z)$ results in a hard singularity at that point. A borderline soft
  singularity has degree of softness zero.

\item A singularity of $w(z)$ at $z_{1}$ is a {\em borderline hard} one if
  it is a hard singularity and a single logarithmic integration of $w(z)$
  results in a soft singularity at that point. A borderline hard
  singularity has degree of hardness zero.

\end{itemize}

\noindent
Finally, the following rules are adopted regarding the superposition of
several singularities as the same point, brought about by the addition of
functions.

\begin{itemize}

\item If there is more than one hard singularity of $w(z)$ superposed at
  $z_{1}$, then the result is a hard singularity and its degree of
  hardness is that of the component singularity with the largest degree of
  hardness.

\item If there is more than one soft singularity of $w(z)$ superposed at
  $z_{1}$, then the result is a soft singularity and its degree of
  softness is that of the component singularity with the smallest degree
  of softness.

\item The superposition of hard singularities and soft singularities of
  $w(z)$ at $z_{1}$ results in a hard singularity at $z_{1}$, with the
  degree of hardness of the component singularity with the largest degree
  of hardness.

\end{itemize}

\noindent
It is not difficult to see that this classification spans all existing
possibilities in so far as the possible types of singularity go. First,
given a point of singularity, either the limit of the function to that
point from within the open unit disk exists or it does not. There is no
third alternative, and therefore every singularity is either soft or
hard. Second, given a soft singularity, either it becomes hard after a
certain finite number of logarithmic differentiations of $w(z)$, or it
does not. Similarly, given a hard singularity, either it becomes soft
after a certain finite number of logarithmic integrations of $w(z)$, or it
does not. In either case there is no third alternative. If the soft or
hard character never changes, then we classify the singularity as
infinitely soft or infinitely hard, as the case may be. Otherwise, we
assign to it a degree of softness or hardness by counting the number $n$
of logarithmic differentiations or logarithmic integrations required to
effect its change of character, and assigning to it the number $n-1$ as
the degree of softness or hardness, as the case may be.

We now recall that there is a set of hard singularities which is already
classified, by means of the concept of the Laurent expansion around an
isolated singular point. If a singularity is isolated in two ways, first
in the sense that there is an open neighborhood around it that contains no
other singularities, and second that it is possible to integrate along a
closed curve around it which is closed in the sense that it does not pass
to another leaf of a Riemann surface when it goes around the point, then
one may write a convergent Laurent expansion for the function around that
point. This leads to the concepts of poles of finite orders and of
essential singularities. In particular, it implies that any analytic
function that has a pole of finite order at the point $z_{1}$ can be
written around that point as the sum of a function which is analytic at
that point and a finite linear combination of the singularities

\begin{displaymath}
  \frac{1}{(z-z_{1})^{n}},
\end{displaymath}

\noindent
for $1\leq n\leq n_{\rm o}$, where $n_{\rm o}$ is the order of the pole.
We can use this set of singularities to illustrate our classification. For
example, if we have an inner analytic function $w(z)$ with a simple pole
$1/(z-z_{1})$ for $z_{1}$ on the unit circle, which is a hard singularity
with degree of hardness $1$, then the logarithmic derivative of $w(z)$ has
a double pole $1/(z-z_{1})^{2}$ at that point, an even harder singularity,
of degree of hardness $2$. Further logarithmic differentiations of $w(z)$
produce progressively harder singularities $1/(z-z_{1})^{n}$, where $n$ is
the degree of hardness of the singularity. We see therefore that multiple
poles fit easily and comfortably into the classification scheme. We may
now proceed to examine this chain of singularities in the other direction,
using logarithmic integration in order to do this.

The logarithmic primitive of the function $w(z)$ mentioned above has a
logarithmic singularity $\ln(z-z_{1})$ at that point, which is the weakest
type of hard singularity in this type of integral-differential chain.
Another logarithmic integration produces a soft singularity such as
$(z-z_{1})\ln(z-z_{1})$, which displays no divergence to infinity. This
establishes therefore that $\ln(z-z_{1})$ is a borderline hard
singularity, with degree of hardness $0$. If we now proceed to
logarithmically differentiate the resulting function, we get back the hard
singularity $\ln(z-z_{1})$. This establishes therefore that
$(z-z_{1})\ln(z-z_{1})$ is a borderline soft singularity, with degree of
softness $0$. This illustrates the transitions between hard and soft
singularities, and also justifies our attribution of degrees of hardness
to the multiple poles, as we did above. Further logarithmic integrations
produce progressively softer singularities such as
$(z-z_{1})^{2}\ln(z-z_{1})$, and so on, were we consider only the hardest
or least soft singularity resulting from each operation and ignore regular
terms, leading to the general expression

\begin{displaymath}
  (z-z_{1})^{n+1}\ln(z-z_{1}),
\end{displaymath}

\noindent
where $n$ is the degree of softness. This completes the examination of the
singularities of this particular type of integral-differential chain. Note
that these soft singularities are isolated in the sense that there is an
open neighborhood around each one of them that contains no other
singularities, but not in the sense that one can integrate in closed
curves around them. This is so because the domains of these functions are
in fact Riemann surfaces with infinitely many leaves, and a curve which is
closed in the complex plane is not really closed in the domain of the
function.

Although this chain of singularities exhausts the possibilities so far as
one is limited to integral-differential chains containing isolated hard
singularities of single-valued functions, there are many other possible
chains of singularities, if one starts with hard singularities having
non-trivial Riemann surfaces, for example such as

\begin{displaymath}
  \frac{\ln(z-z_{1})}{(z-z_{1})^{n}},
\end{displaymath}

\noindent
for $n\geq 1$. One might consider also the more general form

\begin{displaymath}
  (z-z_{1})^{n}\ln^{m}(z-z_{1})
\end{displaymath}

\noindent
for the singularities, where $m\geq 0$ and $n$ in any integer, positive or
negative. This generates quite a large set of possible types of
singularity, both soft and hard.

To complete the picture in our exemplification, a simple and widely known
example of an infinitely hard singularity is an essential singularity such
as $\exp[1/(z-z_{1})]$. On the other hand, an infinitely soft singularity
is not such a familiar object. One interesting example will be discussed
in Section~\ref{APPinfsoft} of Appendix~\ref{APPproofs}.

It is important to note that almost all convergent DP Fourier series will
be related to inner analytic functions either with only soft singularities
on the unit circle or with at most borderline hard singularities, which
will therefore all have non-trivial Riemann surfaces as their domains.
Since in our analysis here we are bound within the unit disk, and will at
most consider limits to the unit circle from within that disk, this is not
of much concern to us, because in this case we never go around one of
these singularities in order to change from one leaf of the Riemann
surface to another. The value of the function $w(z)$ within the unit disk
is defined by its value a the origin, and this determines the leaf of each
Riemann surface which is to be used within the disk. We must always
consider that the branching lines of all such branching points at the unit
circle extend outward from the unit circle, towards infinity.

\section{Convergence and Singularities}

Let us examine the effects on the corresponding series and functions of
each one of these two related types of differentiation and integration
operations in turn. First the effect of differentiations and integrations
with respect to $\theta$ acting on the DP Fourier series, and then the
effect of logarithmic differentiations and integrations acting on the
corresponding inner analytic functions $w(z)$.

On the Fourier side of this discussion, it is clear from the structure of
the DP Fourier series that each differentiation with respect to $\theta$
adds a factor of $k$ to the $a_{k}$ coefficients, and thus makes them go
to zero slower, or not at all, as $k\to\infty$. This either reduces the
rate of convergence of the series or makes them outright divergent.
Integration with respect to $\theta$, on the other hand, has the opposite
effect, since it adds to the $a_{k}$ coefficients a factor of $1/k$, and
thus makes them go to zero faster as $k\to\infty$. This always increases
the rate of convergence of the series. On the inner-analytic side of the
discussion, each logarithmic differentiation increases the degree of
hardness of each singularity, while logarithmic integration decreases it.
At the same time, these operations work in the opposite way on the degrees
of softness. No point of singularity over the unit circle ever vanishes or
appears as a result of these operations. Only the degrees of hardness and
softness change.

We are now in a position to use these preliminary facts to analyze the
question of the convergence of the series on the unit circle. While the
derivations and integrations with respect to $\theta$ change the
convergence status of the DP Fourier series on the unit circle, the
related logarithmic derivations and integrations of the corresponding
inner analytic function do not change the analyticity or the set of
singular points of that function. The only thing that these logarithmic
operations do change is the {\em type} of the singularities over the unit
circle, as described by their degrees of hardness or softness. Hence,
there must be a relation between the mode of convergence or lack of
convergence of the DP Fourier series on the unit circle and the nature of
the singularities of the inner analytic function on that circle. It is
quite clear that the rate of convergence and the very existence of
convergence of the DP Fourier series are tied up to the degree of hardness
or softness of the singularities present. The less soft the singularities,
the slower the convergence, leading eventually to hard singularities and
to the total loss of convergence.

Note that the relation between the singularities of $w(z)$ on the unit
circle and the convergence of the DP Fourier series is non-local, because
processes of differentiation or integration will change the degree of
hardness or softness only locally at the singular points, but will affect
the speed of convergence to zero of the coefficients $a_{k}$, which has
its effects on the rate of convergence of the DP Fourier series everywhere
over the unit circle. Since the relation between the hardness or softness
of the singularities and the convergence of the DP Fourier series is
non-local, the rate of convergence or the lack of it will be ruled by the
hardest or least soft singularity or set of singularities found anywhere
over the whole unit circle. We will call these the {\em dominant
  singularities}. Therefore, from now on we will think in terms of the
dominant singularity or set of singularities which exists on that circle.

The problem of establishing a general, complete and exact set of criteria
determining the convergence or arbitrary DP Fourier series with basis on
the set of dominant singularities of the corresponding inner analytic
function is, so far as we can tell, an open one. We will, however, be able
to classify and obtain the convergence criteria for a fairly large class
of DP Fourier series. In order to establish this class of series, consider
the set of all possible integral-differential chains involving the series
$S_{z}$, $S_{v}$, $S_{\rm c}$ and $S_{\rm s}$, and the respective
functions $w(z)$, $f_{\rm c}(\theta)$ and $f_{\rm s}(\theta)$. Let us use
$S_{v}$ to characterize the elements of these chains. Since the operations
of logarithmic differentiation and logarithmic integration always produce
definite and unique results, it is clear the each series $S_{v}$ belongs
to only one of these chains. Let us now select from the set of all
possible integral-differential chains those that satisfy the following two
conditions.

\begin{itemize}

\item There is in the chain a series $S_{v}$ with coefficients $a_{k}$
  that can be bounded in the following way: there exist positive real
  constants $A_{(-)}$ and $A_{(+)}$, a minimum value $k_{m}$ of $k$ and a
  positive real number $0<\varepsilon<1$ such that for $k>k_{m}$

  \begin{displaymath}
    \frac{A_{(-)}}{k}
    \leq
    |a_{k}|
    \leq
    \frac{A_{(+)}}{k^{\varepsilon}}.
  \end{displaymath}

\item The series $S_{v}$ qualified in the previous condition diverges to
  infinity on at least one point of the unit circle.
 
\end{itemize}

\noindent
We will call the integral-differential chains that satisfy these
conditions {\em regular integral-differential chains}. We will adopt as a
shorthand for the first condition the statement that $|a_{k}|$ behaves as
$1/k^{p}$ for large values of $k$, or $|a_{k}|\propto 1/k^{p}$, where
$0<p\leq 1$. What the condition means is that, while the coefficients go
to zero as $k\to\infty$, they do it sufficiently slowly to prevent the
series $S_{v}$ from being absolutely and uniformly convergent. Therefore
the series $S_{v}$ may still converge, but does not converge absolutely.
Because the series satisfies the second condition we know, from the
extended version of Abel's theorem~\cite{AbelTheo}, that the corresponding
inner analytic function has a divergent limit going to infinity, on at
least one point of the unit circle. Therefore, it has at least one hard
singularity on that circle. We see then that the set of dominant
singularities that the inner analytic function has on the unit circle is
necessarily a set of hard singularities. We will denote the series that
satisfies these conditions in any given regular integral-differential
chain by $S_{v,h0}$, and the corresponding inner analytic function by
$w_{h0}(z)$.

The next series in the chain, obtained from this one by logarithmic
integration, which is the same as integration with respect to $\theta$,
has coefficients that behave as $|a_{k}|\propto 1/k^{p}$ with $1<p\leq 2$,
and is therefore absolutely and uniformly convergent everywhere. It
follows therefore that all the singularities of the corresponding inner
analytic function are soft. We will denote the series obtained from
$S_{v,h0}$ in this way by $S_{v,s0}$, and the corresponding inner analytic
function by $w_{s0}(z)$. The set of dominant singularities of $w_{s0}(z)$
is thus seen to be a set of soft singularities. It follows that the
dominant singularities of $w_{h0}(z)$, which are hard and became soft by
means of a single operation of logarithmic integration, constitute a set
of borderline hard singularities. Since we may go back from $w_{s0}(z)$ to
$w_{h0}(z)$ by a single operation of logarithmic differentiation, it also
follows that the set of dominant singularities of $w_{s0}(z)$ is a set of
borderline soft singularities.

If we go further along in either direction of the chain, in the
integration direction all subsequent series $S_{v,sn}$ are also absolutely
and uniformly convergent, to continuous functions that are everywhere
$C^{n}$ and almost everywhere $C^{n+1}$, typically sectionally $C^{n+1}$,
where $n$ is the degree of softness. The corresponding inner analytic
functions $w_{sn}(z)$ have only soft singularities on the unit circle. In
the other direction, the next series in the chain, denoted by $S_{v,h1}$,
obtained from $S_{v,h0}$ by logarithmic differentiation, has coefficients
that behave as $|a_{k}|\propto k^{p}$ with $0\leq p<1$, and is therefore
divergent everywhere, since the coefficients do not go to zero as
$k\to\infty$. The same can be said of all the subsequent elements
$S_{v,hn}$ of the chain in this direction. The corresponding inner
analytic functions $w_{hn}(z)$ have sets of dominant singularities
consisting of hard singularities. We arrive therefore at the following
scheme of convergence diagnostics, for series $S_{v}$ on regular
integral-differential chains, based on the behavior of the dominant
singularities of $w(z)$.

\begin{itemize}

\item If the dominant singularities of $w(z)$ are hard with degree of
  hardness $n\geq 1$, then the series $S_{v}$ is everywhere divergent.

\item If the dominant singularities of $w(z)$ are borderline hard ones,
  then $S_{v}$ {\em may} still be convergent almost everywhere to a
  discontinuous function, but will diverge in one or more points. If the
  number of dominant singularities is finite, then the limiting function
  will be sectionally continuous and differentiable.

\item If the dominant singularities of $w(z)$ are borderline soft ones,
  then $S_{v}$ is everywhere convergent to an everywhere continuous but
  not everywhere differentiable function. If the number of dominant
  singularities is finite, then the limiting function will be sectionally
  differentiable.

\item If the dominant singularities of $w(z)$ are soft with degree of
  softness $n\geq 1$, then the series $S_{v}$ is everywhere convergent to
  a $C^{n}$ function, which is also $C^{n+1}$ almost everywhere. If the
  number of dominant singularities is finite, then the limiting function
  will be sectionally $C^{n+1}$.

\end{itemize}

\noindent
Given the convergence status of $S_{v}$, corresponding conclusions can
then be drawn for the FC pair of DP Fourier series $S_{\rm c}$ and $S_{\rm
  s}$. Observe that in this argument the facts about the convergence of
the DP Fourier series are in fact feeding back into the question of the
convergence of the power series $S_{z}$ at the rim of its maximum disk of
convergence. It follows therefore that this classification can be
understood as a set of statements purely in complex analysis, since it
also states conditions for the convergence or divergence of the Taylor
series $S_{z}$ over the unit circle, in the cases when that circle is the
boundary of its maximum disk of convergence.

Note that in the case in which the dominant singularities are borderline
hard ones there is as yet no certainty of convergence almost everywhere.
Therefore in this respect this alternative must be left open here, and
will be discussed in the next section for some classes of series. It is an
interesting question whether or not there are series $S_{v}$ with
coefficients that behave as $|a_{k}|\propto 1/k^{p}$ with $0<p\leq 1$ and
that do not diverge to infinity anywhere. The common examples all seem to
diverge as we have assumed here. It seems to be difficult to find an
example that has the opposite behavior, but we can offer no proof one way
of the other. Therefore, we must now leave this here as a mere
speculation.

Observe that, since the inner analytic function can be integrated
indefinitely, as many times as necessary, we have here another way, at
least in principle, to get some information about the function that
originated an arbitrarily given DP Fourier series, besides taking limits
of $w(z)$ from within the open unit disk. We may construct the series
$S_{z}$ from the coefficients and, if the corresponding inner analytic
function $w(z)$ is in fact analytic on the open unit disk, with some upper
bound for the hardness of the singularities on the unit circle, then we
may logarithmically integrate it as many times as necessary to reduce all
the singularities over the unit circle to soft singularities. We will then
have corresponding series $S_{\rm c}$ and $S_{\rm s}$ that are absolutely
and uniformly convergent over the unit circle, with continuous functions
as their limits.

The original function is then a multiple derivative with respect to
$\theta$ of one of these continuous functions. Typically we will not be
able to actually take these derivatives at the unit circle, but we will
always be able to take the corresponding derivatives of $w(z)$ within the
open unit disk, arbitrarily close to the unit circle, so that at least we
will have a chance to understand the origin of the problem, which might
give us an insight into the structure of the application, and related
circumstances, that generated the given DP Fourier series. This procedure
will fail only if the inner analytic function turns out to have an
infinite number of singularities on the unit circle, with degrees of
hardness without an upper bound, or an infinitely hard singularity, such
as an essential singularity, anywhere on the unit circle.

\section{Monotonic Series}

The main remaining problem in our classification of convergence modes in
terms of the dominant singularities is to show that if those singularities
are borderline hard ones then the $S_{v}$ series is still convergent
almost everywhere. While we are not able to prove this in general, even
within the realm of the regular integral-differential chains, there are
some rather large classes of DP Fourier series for which it is possible to
prove the convergence almost everywhere. These series satisfy all the
conditions imposed in the previous section, which means that they belong
to regular integral-differential chains. We will also be able to establish
the existence, location and character of the dominant singularities of the
corresponding inner analytic functions $w(z)$.

We will call these series {\em monotonic series}, which refers to the fact
that they have coefficients $a_{k}$ that behave monotonically with $k$.
Other series can also be built from these monotonic series by means of
finite linear combinations, which will share their properties regarding
convergence and dominant singularities, but which are not themselves
monotonic. We will call these {\em extended monotonic series}. The method
we will use to establish proof of convergence will lead to the concept of
{\em singularity factorization}, which we will later generalize. This is a
method for evaluating the series which is algorithmically useful, and can
be used safely in very general circumstances.

\subsection{A Monotonicity Test}\label{SSECmonotest}

Within the set of all series $S_{v}$ which are very weakly convergent,
there is a subset of series that all converge almost everywhere, as we
will now show. This is the class of series $S_{v}$ that have coefficients
$a_{k}$ that converge monotonically to zero. This monotonicity of the
coefficients can, therefore, be used as a convergence test. No statement
at all has to be made about the speed of convergence to zero, but we may
as well focus our attention on series of type $S_{v,h0}$, with
coefficients behaving as $|a_{k}|\propto 1/k^{p}$ with $0<p\leq 1$, which
are those for which this analysis is most useful, since in this case the
sum

\begin{displaymath}
  \sum_{k=1}^{\infty}
  |a_{k}|
\end{displaymath}

\noindent
diverges to infinity and the series is therefore not absolutely or
uniformly convergent. For simplicity we will take the case in which
$a_{k}>0$ for all $k$, but the argument can be easily generalized to
several other cases, as will be discussed later on. Let us consider then
the series

\begin{displaymath}
  S_{v}
  =
  \sum_{k=1}^{\infty}
  a_{k}v^{k},
\end{displaymath}

\noindent
where $a_{k}>0$ for all $k$, $a_{k+1}\leq a_{k}$ for all $k$ and

\begin{displaymath}
  \lim_{k\to\infty}a_{k}
  =
  0.
\end{displaymath}

\noindent
Such a series can be shown to converge for $\theta\neq 0$ by the use of
the Dirichlet test. For $\theta=0$ it diverges to positive infinity along
the real axis, and therefore satisfies the hypotheses defining a series of
type $S_{v,h0}$ in a regular integral-differential chain. It is not too
difficult to show that, for any non-zero value of $\theta$ in the periodic
interval, the Dirichlet partial sums

\begin{displaymath}
  D_{N}
  =
  \sum_{k=1}^{N}
  v^{k},
\end{displaymath}

\noindent
are contained within a closed disk of radius

\begin{displaymath}
  R
  =
  \frac{1}{2|\sin(\theta/2)|},
\end{displaymath}

\noindent
centered at the point

\begin{displaymath}
  C
  =
  \frac{1}{2}
  +
  \ii\,
  \frac{1}{2\tan(\theta/2)}
\end{displaymath}

\noindent
of the complex plane, for all values of $N$. The absolute values $|D_{N}|$
of the Dirichlet partial sums are therefore bounded by a constant for all
$N$, and since the coefficients go monotonically to zero, the Dirichlet
test applies and the series $S_{v}$ is convergent.

We will, however, demonstrate the convergence in another way, which we
believe to be more fruitful, and more directly related to our problem
here. What we will do is to construct another expression, involving
another DP Fourier series, that converges to the same function and, unlike
$S_{v}$, does so absolutely and uniformly almost everywhere. Consider then
the following algebraic passage-work for our series, for $\theta\neq 0$,
which is equivalent to $v\neq 1$. We start by multiplying and dividing the
series by the factor $(v-1)$, and distributing the one in the numerator,

\noindent
\begin{eqnarray*}
  S_{v}
  & = &
  \frac{v-1}{v-1}
  \sum_{k=1}^{\infty}
  a_{k}v^{k}
  \\
  & = &
  \frac{1}{v-1}
  \left[
    \,
    \sum_{k=1}^{\infty}
    a_{k}v^{k+1}
    -
    \sum_{k=1}^{\infty}
    a_{k}v^{k}
  \right].
\end{eqnarray*}

\noindent
We now change the index of the first series, in order to be able to join
the two resulting series,

\noindent
\begin{eqnarray*}
  S_{v}
  & = &
  \frac{1}{v-1}
  \left[
    \,
    \sum_{k=2}^{\infty}
    a_{k-1}v^{k}
    -
    \sum_{k=1}^{\infty}
    a_{k}v^{k}
  \right]
  \\
  & = &
  \frac{1}{v-1}
  \left[
    -
    a_{1}v
    +
    \sum_{k=2}^{\infty}
    (a_{k-1}-a_{k})v^{k}
  \right].
\end{eqnarray*}

\noindent
If we now define the new coefficients $b_{1}=-a_{1}$ and
$b_{k}=a_{k-1}-a_{k}$ for $k>1$, we have a new series $C_{v}$, which we
name the {\em center series} of $S_{v}$, so that we have

\noindent
\begin{eqnarray*}
  S_{v}
  & = &
  \frac{1}{v-1}\,
  C_{v},
  \\
  C_{v}
  & = &
  \sum_{k=1}^{\infty}
  b_{k}v^{k}.
\end{eqnarray*}

\noindent
The name we chose for this series comes from the fact that it describes
the relatively small drift in the complex plane of the instantaneous
center of rotation of the convergence process of the $S_{v}$ series.
Specially for $a_{k}\propto 1/k^{p}$ with the values of $p$ closer to
zero, the convergence process of the $S_{v}$ series proceeds in a long,
slow spiral around the limiting point, while that of the $C_{v}$ series
goes more or less directly to it. Let us now show that $C_{v}$ is
absolutely and thus uniformly convergent. We simply consider the series
$\overline{C}_{v}$ of the absolute values of the terms of $C_{v}$. Since
$|v|=1$, we get

\noindent
\begin{eqnarray*}
  \overline{C}_{v}
  & = &
  \sum_{k=1}^{\infty}
  |b_{k}|
  \\
  & = &
  a_{1}
  +
  \sum_{k=2}^{\infty}
  |a_{k-1}-a_{k}|.
\end{eqnarray*}

\noindent
Since the coefficients $a_{k}$ decrease monotonically to zero, we have
$|a_{k-1}-a_{k}|=a_{k-1}-a_{k}$ and thus we may write

\noindent
\begin{eqnarray*}
  \overline{C}_{v}
  & = &
  a_{1}
  +
  \sum_{k=2}^{\infty}
  a_{k-1}
  -
  \sum_{k=2}^{\infty}
  a_{k}
  \\
  & = &
  a_{1}
  +
  \sum_{k=1}^{\infty}
  a_{k}
  -
  \sum_{k=2}^{\infty}
  a_{k}
  \\
  & = &
  2a_{1}.
\end{eqnarray*}

\noindent
Since $a_{1}$ is some finite real number, this establishes an upper bound
for $\overline{C}_{v}$. Since this series is a sum of positive terms and
thus monotonically increasing, this implies that it converges and
therefore that $C_{v}$ is absolutely convergent. Since the bound and
therefore the criterion of convergence are independent of $\theta$, the
convergence is also uniform. Therefore, for $v\neq 1$ we may evaluate
$S_{v}$ by first evaluating $C_{v}$ and then multiplying the result by the
simple pole $1/(v-1)$.

This proves not only the absolute and uniform convergence of the series
$C_{v}$ everywhere on the unit circle, it also implies the convergence of
the original series $S_{v}$ at all points except the special point $v=1$,
which corresponds to $\theta=0$. Obviously the $S_{v}$ series does not
converge absolutely or uniformly, but it does converge at all points
except $v=1$. We see therefore that the monotonicity of the coefficients
can be used as a test for simple point-wise convergence of the $S_{v}$
series almost everywhere on the unit circle. The special point, were the
$S_{v}$ series diverges, is easily identified as the point where $w(z)$
has a borderline hard singularity, since a single logarithmic integration
of $w(z)$ will necessarily result in a soft singularity there.

If we extend the series $C_{v}$ to a full complex power series $C_{z}$ on
the unit disk, just as we did before in the case of $S_{v}$ and $S_{z}$,
we immediately see that

\begin{displaymath}
  C_{z}
  =
  (z-1)S_{z}.
\end{displaymath}

\noindent
Note that $C_{z}$ satisfies all the necessary conditions for convergence
to an inner analytic function. Since $S_{z}$ converges to an analytic
function on the open unit disk and $(z-1)$ is analytic on the whole
complex plane, it follows that $C_{z}$ also converges to an analytic
function on that disk. In addition to this, since $(z-1)$ reduces to a
real function on the real line and $S_{z}$ reduces to a real function on
the real interval $(-1,1)$, so does $C_{z}$. Finally, it is clear that
since $S_{z}$ is zero at $z=0$, so is $C_{z}$. Therefore $C_{z}$ converges
to an inner analytic function $\gamma(z)$, so that we have

\begin{displaymath}
  \gamma(z)
  =
  (z-1)w(z).
\end{displaymath}

\noindent
Since $C_{z}$ is absolutely and uniformly convergent on the unit circle,
the function $\gamma(z)$ can have only soft singularities on that circle.
However, $S_{z}$ diverges to infinity at one point on the unit circle and,
due to the extended version of Abel's theorem~\cite{AbelTheo}, it follows
that $w(z)$ must have a hard singularity at that point. We see therefore
that in this case the multiplication by $(z-1)$ has the same effect of a
logarithmic integration. While the function $w(z)$ has a single borderline
hard singularity at $z=1$, the function $\gamma(z)$ has a borderline soft
singularity at that point.

We may conclude here that this whole class of series satisfies our
hypotheses defining a series of type $S_{v,h0}$ in a regular
integral-differential chain, as well as that all of them converge almost
everywhere. Besides, all the series in this whole class have a single
dominant borderline hard singularity located at $z=1$. Note that even if a
monotonic series has coefficients that go to zero as a power in ways other
than $1/k^{p}$ with $0<p\leq 1$, the construction of $C_{z}$ still
applies. In this case $w(z)$ will still have a dominant singularity at
$z=1$, although no longer a borderline hard one, but a soft one instead.

\subsection{Monotonicity Extensions}

The argument given in the last subsection can be generalized in several
ways, and the generalizations constitute proof of the convergence of wider
classes of series of type $S_{v,h0}$. For example, it can be trivially
generalized to series with negative $a_{k}$ for all $k$, converging to
zero from below. It is also trivial that it can be generalized to series
which include a non-monotonic initial part, up to some minimum value
$k_{m}$ of $k$. In a less trivial, but still simple way, the argument can
be generalized to series with only odd-$k$ terms, such as the square wave,
or with only even-$k$ terms, such as the two-cycle sawtooth wave. We can
prove this in a simple way by reducing these cases to the previous one. If
we have a complex power series given by

\begin{displaymath}
  S_{z}
  =
  \sum_{j=1}^{\infty}
  a_{k}z^{k},
\end{displaymath}

\noindent
where $k=2j$, we may simply define new coefficients $a'_{j}=a_{k}$ and a
new complex variable $z'=z^{2}$, in terms of which the series can now be
written as

\begin{displaymath}
  S_{z}
  =
  \sum_{j=1}^{\infty}
  a'_{j}z^{\prime j},
\end{displaymath}

\noindent
thus reducing it to the previous form, in terms of the variable $z'$. So
long as the non-zero coefficients $a_{k}$ tend monotonically to zero as
$k\to\infty$, we have that the coefficients $a'_{j}$ tend monotonically to
zero as $j\to\infty$, and the previous result applies. The same is true if
we have a complex power series given by

\begin{displaymath}
  S_{z}
  =
  \sum_{j=0}^{\infty}
  a_{k}z^{k},
\end{displaymath}

\noindent
where $k=2j+1$, since we may still define new coefficients $a'_{j}=a_{k}$
and the new complex variable $z'=z^{2}$, in terms of which the series can
now be written as

\begin{displaymath}
  S_{z}
  =
  z
  \sum_{j=0}^{\infty}
  a'_{j}z^{\prime j},
\end{displaymath}

\noindent
so that once more, so long as the non-zero coefficients tend monotonically
to zero, the previous result applies. We say that series such as these
have non-zero coefficients with a constant step $2$. The only important
difference that comes up here is that the special points on the unit
circle are now defined by $z'=1$, which means that $z^{2}=1$ and hence
that $z=\pm 1$. Therefore, in these cases one gets two special points on
the unit circle, instead of one, namely $\theta=0$ and $\theta=\pm\pi$, at
which we have dominant singularities, which will be borderline hard
singularities so long as $a_{j}\propto 1/j^{p}$ with $0<p\leq 1$.

In addition to this, series with step $1$ and coefficients that have
alternating signs, such as $a_{k}=(-1)^{k}b_{k}$ with monotonic $b_{k}$,
which are therefore not monotonic series, can be separated into two
sub-series with step $2$, one with odd $k$ and the other with even $k$,
and since the coefficients of these two sub-series are monotonic, then the
result holds for each one of the two series, and hence for their sum. In
this case we will have two dominant singularities in each one of the
components series, located at $z=\pm 1$. However, sometimes the two
singularities at $z=1$, one in each component series, may cancel off and
the original series may have a single dominant singularity located at
$z=-1$. One can see this by means of a simple transformation of variables,

\begin{displaymath}
  (-1)^{k}b_{k}z^{k}
  =
  b_{k}z^{\prime k},
\end{displaymath}

\noindent
where $z'=-z$, so that $z'=1$ implies $z=-1$. Series with step $2$ and
coefficients that have alternating signs, such as $a_{k}=(-1)^{j}b_{k}$
with $k=2j$ or $k=2j+1$ and monotonic $b_{k}$, which are also not
monotonic themselves, can be separated into two sub-series with step $4$,
and since the coefficients of these two sub-series are monotonic, then the
result holds for each one of the two series, and hence for their sum. In
this case we will have four dominant singularities in each one of the
components series, located at $z=\pm 1$ and $z=\pm\ii$. However, sometimes
the singularities at $z=\pm 1$ of the two component series may cancel off
and the original series may have only two dominant singularity located at
$z=\pm\ii$. One can see this by means of another simple transformation of
variables,

\begin{displaymath}
  (-1)^{j}b_{k}z^{2j}
  =
  b_{k}z^{\prime j},
\end{displaymath}

\noindent
where $z'=-z^{2}$, so that $z'=1$ implies $z^{2}=-1$ and hence $z=\pm\ii$.
In fact, the result can be generalized to series with non-zero terms only
at some arbitrary regular interval $\Delta k$, that is, having non-zero
terms with some constant step other than $2$. If we have a complex power
series given by

\begin{displaymath}
  S_{z}
  =
  \sum_{j=0}^{\infty}
  a_{k}z^{k},
\end{displaymath}

\noindent
where $k=k_{0}+pj$, for some strictly positive integer $k_{0}$ and where
the step $p$ is another strictly positive integer, we may simply define
new coefficients $a'_{j}=a_{k}$ and a new complex variable $z'=z^{p}$, in
terms of which the series can now be written as

\begin{displaymath}
  S_{z}
  =
  z^{k_{0}}
  \sum_{j=0}^{\infty}
  a'_{j}z^{\prime j},
\end{displaymath}

\noindent
thus reducing it to the previous form, in terms of the variable $z'$. So
long as the non-zero coefficients $a_{k}$ tend monotonically to zero as
$k\to\infty$, we have that the coefficients $a'_{j}$ tend monotonically to
zero as $j\to\infty$, and the previous result applies. In this case the
special points over the unit circle are given by $z^{p}=1$, and there are,
therefore, $p$ such points, including $z=1$, uniformly distributed along
the circle. If combined with alternating signs, these series have special
points given by $z^{p}=-1$, and once again there are $p$ such points
uniformly distributed along the circle. Note that the number of dominant
singularities on the unit circle increases with the step $p$, and that
they are homogeneously distributed along that circle.

Finally, one may consider building finite superpositions of the series in
all the previous cases discussed so far. Since each component series
converges almost everywhere and has a finite number of dominant
singularities, these superpositions of series will all converge, and will
all have a finite number of dominant singularities on the unit disk. Since
all the component series are of type $S_{v,h0}$, and hence have dominant
singularities which are borderline hard ones, the dominant singularities
of the superpositions will always be at most borderline hard ones. If the
dominant singularities are in fact borderline hard ones, then we will call
these series {\em extended monotonic series}. Each one of these series
generates a different regular integral-differential chain, and this
defines a rather large set of series that can be classified according to
our scheme, relating their mode of convergence and the dominant
singularities on the unit circle.

\section{Factorization of Singularities}

One of the interesting facts that follow from the analysis in the previous
paper~\cite{FTotCPI} is that, since the limiting function of a DP Fourier
series is always given by the limit of the corresponding inner analytic
function from within the open unit disk to the unit circle, in any open
subset of that circle where $w(z)$ is analytic it is also $C^{\infty}$
along $\theta$. Therefore a DP Fourier series that converges in a
piecewise fashion between two consecutive singularities of $w(z)$ does so
to a piecewise section of a $C^{\infty}$ function. This means that it
should be possible to recover the $C^{\infty}$ function involved in each
section, and also that it should be possible to represent them by series
that converge at a faster rate and can thus be differentiated at least a
few times. In this section we will show how one can accomplish the latter
goal.

We start with a simple case, which in fact we have already demonstrated
completely in the previous section. The proof of convergence of DP Fourier
series with monotonic coefficients described in
Subsection~\ref{SSECmonotest} can be understood as a process of
factorization of the singularity of the inner analytic function $w(z)$.
Interpreted in terms of $S_{z}$ we may write the relation between that
series and the corresponding center series $C_{z}$ as

\noindent
\begin{eqnarray*}
  S_{z}
  & = &
  \frac{1}{z-1}\,
  C_{z},
  \\
  C_{z}
  & = &
  \sum_{k=1}^{\infty}
  b_{k}z^{k}.
\end{eqnarray*}

\noindent
As we have shown before, since $S_{z}$ converges to an inner analytic
function $w(z)$, so does $C_{z}$, and hence we have

\begin{displaymath}
  w(z)
  =
  \frac{1}{z-1}\,
  \gamma(z),
\end{displaymath}

\noindent
where $C_{z}$ converges to $\gamma(z)$. What was done here is to factor
out of $S_{z}$ a simple pole at the point $z=1$. Hence the original series
$S_{v}$, which is not absolutely or uniformly convergent and is associated
to an inner analytic function that has a borderline hard singularity at
$z=1$, is translated into a series $C_{v}$ which is absolutely and
uniformly convergent, and that is associated to an inner analytic function
that has a borderline soft singularity at that point. Note that the
$S_{z}\to C_{z}$ transformation does not change the maximum disk of
convergence or the location of any singularities. Just like logarithmic
integration, it just softens the existing singularities.

\subsection{General Singularity Factorization}

If we think about our general scheme of classification of singularities
and modes of convergence, we can see that so long as that scheme holds
this process of factoring out singularities should always work, regardless
of any hypothesis about the coefficients, such as that they be monotonic.
So long as the coefficients of the original DP trigonometric series lead
to the construction of an inner analytic function, and so long as that
inner analytic function has at most a finite set of dominant singularities
over the unit circle, which are not infinitely hard ones such as essential
singularities, it should be possible to do this, and hence produce another
related series in which the dominant singularities are softened.

Here is how this procedure works. Given a certain DP Fourier series with
coefficients $a_{k}$, we construct the series $S_{z}$ and thus the inner
analytic function $w(z)$ and determine the set of dominant singularities
that it has on the unit circle, which we assume are $N$ in
number. Independently of the degree of hardness or softness of these
singularities, we now introduce simple poles at each dominant singularity,

\noindent
\begin{eqnarray*}
  S_{z}
  & = &
  \frac
  {(z-z_{1})\ldots(z-z_{N})}
  {(z-z_{1})\ldots(z-z_{N})}\,
  S_{z}
  \\
  & = &
  \frac{1}{(z-z_{1})\ldots(z-z_{N})}\,
  C_{z},
\end{eqnarray*}

\noindent
where the new series is defined by

\begin{displaymath}
  C_{z}
  =
  P_{N}(z)
  S_{z},
\end{displaymath}

\noindent
where $P_{N}(z)$ is a polynomial of the order indicated,

\begin{displaymath}
  P_{N}(z)
  =
  (z-z_{1})\ldots(z-z_{N}).
\end{displaymath}

\noindent
This polynomial and the original series $S_{z}$ can then be manipulated
algebraically in order to produce an explicit expression for the new
series $C_{z}$, which we will still call the center series of $S_{z}$. Let
us show that $C_{z}$ converges to an inner analytic function $\gamma(z)$.
Since $S_{z}$ converges to an analytic function $w(z)$ on the open unit
disk, it is at once apparent that $C_{z}$ also converges to an analytic
function on that disk, since $P_{N}(z)$ is a polynomial and hence an
analytic function over the whole complex plane. Also, since $S_{z}=0$ at
the point $z=0$, it follows at once that $C_{z}=0$ on that same point.

Let us now recall that $S_{z}$ is a power series generated by a FC pair of
DP Fourier series, and therefore that its real and imaginary parts have
definite parities. Therefore, the inner analytic function $w(z)$ that it
converges to also has real and imaginary parts with definite parities. As
we showed before, its real part is even on $\theta$, and its imaginary
part is odd on $\theta$. Therefore, the singularities of the function
$w(z)$ must come in pairs, unless they are located at $\theta=0$ or
$\theta=\pm\pi$. This means that, if there is a singularity at a point
$z_{1}$ on the unit circle away from the real axis, then there is an
essentially identical one at $z_{1}^{*}$,

\noindent
\begin{eqnarray*}
  z_{1}
  & = &
  \cos(\theta_{1})
  +
  \ii
  \sin(\theta_{1})
  \;\;\;\Rightarrow
  \\
  z_{1}^{*}
  & = &
  \cos(\theta_{1})
  -
  \ii
  \sin(\theta_{1})
  \\
  & = &
  \cos(-\theta_{1})
  +
  \ii
  \sin(-\theta_{1}),
\end{eqnarray*}

\noindent
possibly with the overall sign reversed. It follows that, if we want to
factor out the singularities on both points, we must choose the factors
that constitute $P_{N}(z)$ in pairs of factors at complex-conjugate
points, except possibly for a couple of points over the real axis.
Assuming for example that $z_{2}$ is real, we have to use something like

\begin{displaymath}
  P_{N}(z)
  =
  (z-z_{1})
  (z-z_{1}^{*})
  \times
  (z-z_{2})
  \times
  (z-z_{3})
  (z-z_{3}^{*})
  \times
  \ldots\;.
\end{displaymath}

\noindent
If we restrict the polynomial to the real axis we get

\begin{displaymath}
  P_{N}(x)
  =
  (x-z_{1})
  (x-z_{1}^{*})
  \times
  (x-z_{2})
  \times
  (x-z_{3})
  (x-z_{3}^{*})
  \times
  \ldots\;.
\end{displaymath}

\noindent
If we now take the complex conjugate of $P_{N}(x)$ we see that in fact
nothing changes,

\begin{displaymath}
  P_{N}^{*}(x)
  =
  (x-z_{1}^{*})
  (x-z_{1})
  \times
  (x-z_{2})
  \times
  (x-z_{3}^{*})
  (x-z_{3})
  \times
  \ldots\;.
\end{displaymath}

\noindent
Since we thus conclude that $P_{N}^{*}(x)=P_{N}(x)$, it follows that
$P_{N}(x)$ is a real polynomial over the real axis. Since the series
$S_{z}$ converges to an inner analytic function $w(z)$, which also reduces
to a real function on the interval $(-1,1)$ of the real axis, it follows
that $\gamma(z)$ reduces to a real function on the interval $(-1,1)$ of
the real axis as well. This establishes that the function $\gamma(z)$ has
all the required properties and is therefore an inner analytic function.

Assuming that the series $S_{z}$ is convergent, and thus that $w(z)$ has
at most borderline hard singularities, this new series $C_{z}$ generates a
new inner analytic function $\gamma(z)$ that has only soft singularities
and hence $C_{z}$ converges absolutely and uniformly to a continuous
function. From that function and the explicit poles we can then
reconstruct the original function, in a piecewise fashion between pairs of
adjacent dominant singularities. If the series $S_{z}$ was already
absolutely and uniformly convergent, the new center series $C_{z}$ will
allow one to take one more derivative, compared with the situation
regarding $S_{z}$.

Note that between two adjacent singularities the function $w(z)$ is
analytic over sections of the unit circle, and hence piecewise
$C^{\infty}$, so that this process can be taken, in principle, as far as
one wishes, by the iteration of this procedure. In order to do this one
has to re-examine the set of singularities of the resulting function
$\gamma(z)$ because, with the softening of the dominant singularities,
there may be now more singularities just as soft as those in the first set
became. This will generate a new set of dominant singularities, and
assuming that this set is also finite in number, one may iterate the
procedure. The result will be a series that not only is absolutely and
uniformly convergent, but is also one that can be differentiated one more
time and still result in an equally convergent series.

Note also that we may as well start the process with a series that if
flatly divergent, and that there is nothing to prevent us from recovering
from it the original function, from which the coefficients were obtained.
This can always be done, at least in principle, if the inner analytic
function has at most a finite number of isolated singularities on the unit
circle, each one with a finite degree of hardness. On the other hand, it
cannot be done if there is an infinite number of dominant singularities,
or if any individual singularity is an infinitely hard one, such as an
essential singularity.

Arguably the most difficult step in this process is the determination of
the singularities of the inner analytic function from the series. If there
is enough information about the real function that originated it, or about
the circumstances of the application involved, it may be possible to guess
at the set of singularities. Otherwise, this information has to be
obtained from the structure of the series itself.

However, this may not be such a grave difficulty as it might appear at
first, since a given singularity structure characterizes a whole class of
series, not a single series. For example, all series which have monotonic
coefficients with step $1$ have a single dominant singularity at $z=1$.
Series which have monotonic coefficients with step $1$ and a factor of
$(-1)^{k}$ added to the coefficients have a single dominant singularity at
$z=-1$. Series with monotonic coefficients with step $2$ have two dominant
singularities, at $z=1$ and at $z=-1$. Series with monotonic coefficients
with step $2$ and a factor of $(-1)^{j}$ added to the coefficients have
two dominant singularities, at $z=\ii$ and at $z=-\ii$, and so on. In
Appendix~\ref{APPcenter} we will give several simple examples of the
construction of center series.

Observe also that this whole procedure is safe in the sense that if one
guesses erroneously at the singularities, the worst that can happen is
that no improvement in convergence is obtained. Besides, at least in
principle the factorization process can be considered in reverse, in the
sense that one may analyze the structure of the series $S_{z}$ in order to
discover what set of factors would do the trick of resulting in a series
$C_{z}$ with coefficients that go to zero faster than the original ones.
If this problem is solvable, it in fact {\em determines} the location of
the dominant singularities of the inner analytic function over the unit
circle, by what turns out to have the nature of a purely algebraic method,
leading to a polynomial $P_{N}(z)$ that implements the softening of the
dominant singularities.

\section{Conclusions}

The sometimes complicated questions of convergence of Fourier series can
be mapped onto the convergence of Taylor series of analytic functions.
The modes of convergence of DP Fourier series can be classified according
to the singularity structure of the corresponding inner analytic
functions. The extreme cases of strong convergence and strong divergence
are easily identified and classified, as was shown in the previous
paper~\cite{FTotCPI}. Simple tests can be used to identify these cases.

Weakly convergent DP Fourier series present a much more delicate and
difficult problem. It can be shown that all these cases are translated, in
the complex formalism, into the behavior of inner analytic functions and
their Taylor series at the rim of the maximum disk of convergence of these
series, which in this case is the unit circle. The treatment of this case
required the introduction of a classification scheme for singularities of
inner analytic functions. These were classified as either soft or hard,
depending on the behavior of the inner analytic functions near them, and
subsequently by integer indices giving the degrees of either softness or
hardness of the singularities. Of particular importance are the degrees
which we named as borderline soft and borderline hard.

This classification scheme led to the concept of integral-differential
chains of inner analytic functions, which were needed in order to relate
the classification of singularities with a corresponding classification of
modes of convergence of the series associated to the functions. Definite
results were obtained only for a certain subset of all possible such
chains, consisting of regular chains in which the series $S_{v,h0}$ is an
extended monotonic series. This defines a certain class of series and
corresponding functions. With the limitation that we must stay within this
class, it was shown that the existence and level of convergence of DP
Fourier series is ruled by the nature of the dominant singularities of the
inner analytic functions which are located at the rim of their maximum
disk of convergence, which is the unit circle.

As a result of this classification, by constructing the inner analytic
function of a given DP Fourier series in this class, one can determine the
convergence of the series via the examination of the singularities of that
function. At a very basic level, one can determine whether the series is
strongly divergent or strongly convergent by simply determining the
position of possible singularities of the inner analytic function. This
leads to the three-pronged basic decision process: if there is a
singularity of any type within the open unit disk, then the $S_{v}$ series
is divergent everywhere; if there are no singularities within the closed
unit disk, then the $S_{v}$ series converges everywhere to a $C^{\infty}$
function; if there are one or more singularities on the unit circle, but
none within the open unit disk, then the convergence of the $S_{v}$ series
is determined by the degree of hardness or of softness of the dominant
singularities on the unit circle.

The last alternative leads to another three-pronged decision process, this
time based on the type of the dominant singularities of the inner analytic
function found on the unit circle. According to the underlying structure
that was uncovered, at this finer level the decision structure leads to
the following basic alternatives: if the dominant singularities are soft
singularities, then the $S_{v}$ series converges absolutely and uniformly
everywhere; if they are borderline hard singularities, then the $S_{v}$
series converges point-wise almost everywhere, but does not converge
absolutely; if they are hard singularities with degree of hardness $1$ or
greater, then the series $S_{v}$ diverges everywhere.

This classification scheme for all the modes of convergence, and of the
corresponding degrees of softness or hardness of the dominant
singularities, is illustrated in Table~\ref{TABclass}, which gives also
some additional information. By and large, as the singularities become
softer the series become more convergent over the unit circle, and
converge to smoother functions. Given the convergence mode of the $S_{v}$
series, one can then derive the corresponding mode for the DP Fourier
series, which are also included in the table.

\begin{table}[htp]
  \centering
  \begin{tabular}{||c|c|c|c|c||}
    \hline
    \hline
    Dominant				&
    Convergence				&
    Behavior of				&
    Convergence				&
    Character
    \\
    Singularities			&
    of $S_{v}$				&
    Coefficients			&
    of $S_{\rm c}$ and $S_{\rm s}$	&
    of $f(\theta)$
    \\
    \hline
    \hline
    $n$-hard				&
    divergent				&
    $|a_{k}|\propto k^{p}$		&
    divergent				&
    currently
    \\
    $n\geq 2$				&
    ew					&
    $n-1\leq p<n$			&
    aew					&
    unknown
    \\
    \hline
    $1$-hard				&
    divergent				&
    $|a_{k}|\propto k^{p}$		&
    divergent				&
    $\delta$-``function''
    \\
    \mbox{}				&
    ew					&
    $0\leq p<1$				&
    aew					&
    for $p=0$
    \\
    \hline
    borderline				&
    point-wise				&
    $|a_{k}|\propto 1/k^{p}$		&
    point-wise				&
    cont aew
    \\
    hard				&
    aew					&
    $0<p\leq 1$				&
    aew					&
    diff aew
    \\
    \hline
    borderline				&
    absolute and			&
    $|a_{k}|\propto 1/k^{p}$		&
    absolute and			&
    cont ew
    \\
    soft				&
    uniform ew				&
    $1<p\leq 2$				&
    uniform ew				&
    diff aew
    \\
    \hline
    $1$-soft				&
    absolute and			&
    $|a_{k}|\propto 1/k^{p}$		&
    absolute and			&
    diff ew
    \\
    \mbox{}				&
    uniform ew				&
    $2<p\leq 3$				&
    uniform ew				&
    $C^{2}$ aew
    \\
    \hline
    $n$-soft				&
    absolute and			&
    $|a_{k}|\propto 1/k^{p}$		&
    absolute and			&
    $C^{n}$ ew
    \\
    $n\geq 2$				&
    uniform ew				&
    $n+1<p\leq n+2$			&
    uniform ew				&
    $C^{n+1}$ aew
    \\
    \hline
    \hline
  \end{tabular}
  \caption{A table showing the proposed classification of modes of
    convergence, singularity structure and limiting function properties,
    within the class of regular integral-differential chains in which the
    series $S_{v,h0}$ is an extended monotonic series. The abbreviation
    ``ew'' stands for ``everywhere'' and ``aew'' for ``almost
    everywhere''. The abbreviations ``cont'' and ``diff'' stand 
    respectively for ``continuous'' and ``differentiable''. Integration
    goes downward through the lines, differentiation goes
    upward.}\label{TABclass}
\end{table}

In terms of the analytic character of the limiting functions, at the most
basic level strongly convergent DP Fourier series converge to restrictions
to the whole unit circle of complex $C^{\infty}$ analytic functions, while
weakly convergent DP Fourier series converge to globally $C^{n}$ functions
which are also sectionally $C^{n+1}$, where $n$ is the degree of softness
of the dominant singularities on the unit circle. In the case in which the
dominant singularities are borderline hard, the series converge to
sectionally continuous and differentiable functions, which however are not
globally continuous. In addition to this, in the case in which the
dominant singularities are simple poles one may have representations of
singular objects such as the Dirac delta ``function'', as was shown in the
previous paper~\cite{FTotCPI}. However, this last alternative has not yet
been explored in much detail.

Moreover, we presented the process of singularity factorization, through
which, given an arbitrary DP Fourier series, which can even be divergent
almost everywhere, one can construct from it other expressions involving
trigonometric series, that converge to the function that gave origin to
the given DP Fourier series. This works by the construction of a new
complex power series from the Taylor series $S_{z}$ of the corresponding
inner analytic function. We call these new series $C_{z}$ the center
series of the series $S_{z}$ that converges to the inner analytic function
associated to the original DP Fourier series. If the original series was
very weakly convergent, then this new series will have much better
convergence characteristics. Even if the original series is divergent one
can still construct expressions involving center series that converge to
the original function, in a piecewise fashion. In this way, a more
practical means of recovery of the original function is provided, if
compared to the explicit determination of the inner analytic function
$w(z)$ in closed form, in order to enable one to take its limit to the
unit circle explicitly.

Several points are left open and represent interesting possibilities for
further development of the subject. One was presented in the previous
paper~\cite{FTotCPI} and consists of the question of whether or not there
are real functions which generate strongly divergent DP Fourier series.
The conjecture is that there are none, in which case taking limits of
inner analytic functions from within the open unit disk would be
established as a process for the generation, almost everywhere on the unit
circle, of all real functions from which it is possible to define the
coefficients of a DP Fourier series. Another interesting question, posed
in this paper, is whether or not there are $S_{v}$ series with
coefficients $a_{k}$ that behave as $|a_{k}|\propto 1/k^{p}$ with $0<p\leq
1$ and that converge everywhere on the unit circle. If there are, it would
be necessary to consider the extension of our classification scheme to
other classes of DP Fourier series.

Since an absolutely and uniformly convergent DP Fourier series usually
converges {\em much} faster than a non-absolutely and non-uniformly
convergent one, doing the $S_{v}\to C_{v}$ transformation can be of
enormous numerical advantage. One verifies that, the slower the
convergence of $S_{v}$, caused by a value of $p$ closer to zero when
$a_{k}\propto 1/k^{p}$, the more advantageous is the use of the series
$C_{v}$. Near the special points the gain is more limited, but it still
exists. For simple well-known series such as the square wave, with
coefficients that go to zero as $1/k$, which is the fastest possible
approach within the $S_{v,h0}$ class, on average over the whole domain,
and for the higher levels of numerical precision required of the results,
the speedup can be as high as $1000$ or more, as we will show
elsewhere~\cite{CSNumerics}.

\section{Acknowledgements}

The author would like to thank his friend and colleague Prof. Carlos
Eugênio Imbassay Carneiro, to whom he is deeply indebted for all his
interest and help, as well as his careful reading of the manuscript and
helpful criticism regarding this work.

\appendix

\section{Appendix: Technical Proofs}\label{APPproofs}

\subsection{Evaluations of Convergence}\label{APPevalconv}

It is not a difficult task to establish the absolute and uniform
convergence of DP Fourier series, or the lack thereof, starting from the
behavior of the coefficients of the series in the limit $k\to\infty$, if
we assume that they behave as inverse powers of $k$ for large values of
$k$. If we have a complex series $S_{v}$ with coefficients $a_{k}$,

\begin{displaymath}
  S_{v}
  =
  \sum_{k=1}^{\infty}a_{k}v^{k},
\end{displaymath}

\noindent
where $|v|=1$, then it is absolutely convergent if and only if the series
$\overline{S}_{v}$ of the absolute values of the coefficients,

\noindent
\begin{eqnarray*}
  \overline{S}_{v}
  & = &
  \sum_{k=1}^{\infty}|a_{k}||v|^{k}
  \\
  & = &
  \sum_{k=1}^{\infty}|a_{k}|,
\end{eqnarray*}

\noindent
converges. One can show that this sum will be finite if, for $k$ above a
certain minimum value $k_{m}$, it holds that

\begin{displaymath}
  |a_{k}|
  \leq
  \frac{A}{k^{1+\varepsilon}},
\end{displaymath}

\begin{figure}[ht]
  \centering
  \fbox{
    \epsfig{file=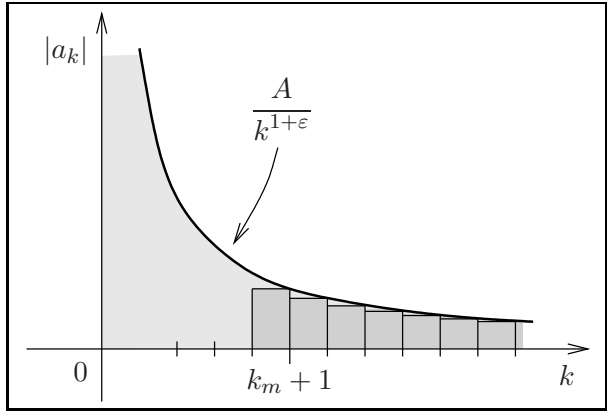,scale=1.0,angle=0}
  }
  \caption{Illustration of the upper bounding of the sum of the
    coefficients $|a_{k}|$ by an integral.}
  \label{APPfig01}
\end{figure}

\noindent
for some positive real constant $A$ and some real constant
$\varepsilon>0$. This is true because the sum of a finite set of initial
terms is necessarily finite, and because in this case we may bound the
remaining infinite sum from above by a convergent asymptotic integral,

\noindent
\begin{eqnarray*}
  \sum_{k=k_{m}+1}^{\infty}|a_{k}|
  & \leq &
  \sum_{k=k_{m}+1}^{\infty}\frac{A}{k^{1+\varepsilon}}
  \\
  & < &
  \int_{k_{m}}^{\infty}dk\,\frac{A}{k^{1+\varepsilon}}
  \\
  & = &
  \frac{-A}{\varepsilon}\,\frac{1}{k^{\varepsilon}}\at{k_{m}}{\infty}
  \\
  & = &
  \frac{A}{\varepsilon}\,\frac{1}{k_{m}^{\varepsilon}},
\end{eqnarray*}

\noindent
as illustrated in Figure~\ref{APPfig01}. In that illustration each
vertical rectangle has base $1$ and height given by $|a_{k}|$, and
therefore area given by $|a_{k}|$. As one can see, the construction is
such that the set of all such rectangles is below the graph of the
function $A/k^{1+\varepsilon}$, and therefore the sum of their areas is
contained within the area under that graph, to the right of $k_{m}$. This
establishes the necessary inequality between the sum and the integral.

So long as $\varepsilon$ is not zero, this establishes an upper bound to a
sum of positive quantities, which is therefore a monotonically increasing
sum. It then follows from the well-known theorem of real analysis that the
sum necessarily converges, and therefore the series $S_{v}$ is absolutely
convergent. The same is then true for the corresponding DP Fourier series.

In addition to this, one can see that the convergence condition does not
depend on $\theta$, since that dependence is only within the complex
variable $v=\exp(\ii\theta)$, and vanishes when we take absolute values.
This implies uniform convergence because, given a strictly positive real
number $\epsilon$, absolute convergence for this value of $\epsilon$
implies convergence for this same value of $\epsilon$, with the same
solution $k(\epsilon)$ for the convergence condition. This makes it clear
that the solution of the convergence condition for $k$ is independent of
position and therefore that the series is also uniformly convergent. Once
more, the same is then true for the corresponding DP Fourier series.

\begin{figure}[ht]
  \centering
  \fbox{
    \epsfig{file=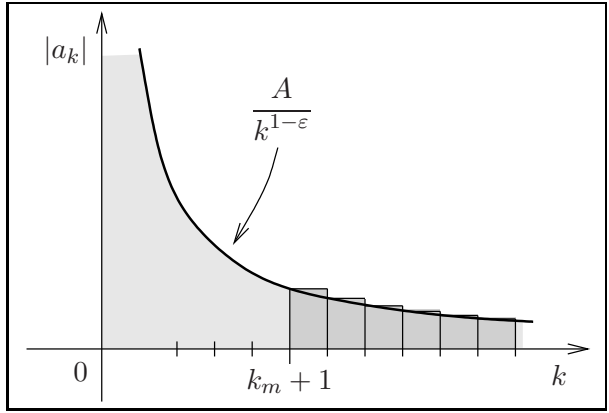,scale=1.0,angle=0}
  }
  \caption{Illustration of the lower bounding of the sum of the
    coefficients $|a_{k}|$ by an integral.}
  \label{APPfig02}
\end{figure}

This establishes a sufficient condition for the absolute and uniform
convergence of DP Fourier series. On the other hand, if we have that, for
$k$ above a certain minimum value $k_{m}$,

\begin{displaymath}
  |a_{k}|
  \geq
  \frac{A}{k^{1-\varepsilon}},
\end{displaymath}

\noindent
with positive real $A$ and real $\varepsilon\geq 0$, then it is possible
to bound the sum $\overline{S}_{v}$ from below by an asymptotic integral
that diverges to positive infinity. This is done in a way similar to the
one used for the establishment of the upper bound, but inverting the
situation so as to keep the area under the graph contained within the
combined areas of the rectangles, as illustrated in Figure~\ref{APPfig02}.
The argument then establishes in this case that, for $\varepsilon>0$

\noindent
\begin{eqnarray*}
  \sum_{k=k_{m}+1}^{\infty}|a_{k}|
  & \geq &
  \sum_{k=k_{m}+1}^{\infty}\frac{A}{k^{1-\varepsilon}}
  \\
  & > &
  \int_{k_{m}+1}^{\infty}dk\,\frac{A}{k^{1-\varepsilon}}
  \\
  & = &
  \frac{A}{\varepsilon}\,k^{\varepsilon}\at{k_{m}+1}{\infty}
  \\
  & = &
  -\,
  \frac{A}{\varepsilon}\,(k_{m}+1)^{\varepsilon}
  +
  \frac{A}{\varepsilon}
  \lim_{k\to\infty}
  k^{\varepsilon},
\end{eqnarray*}

\noindent
and therefore that $\overline{S}_{v}$ diverges to infinity. A similar
calculation can be performed in the case $\varepsilon=0$, leading to
logarithms and yielding the same conclusions. This does not prove or
disprove convergence itself, but it does establish the absence of absolute
convergence. It also shows that, so long as $|a_{k}|$ behaves as a power
of $k$ for large $k$, the previous condition is both sufficient and
necessary for absolute convergence.

\subsection{Infinitely Soft Singularities}\label{APPinfsoft}

Consider the following power series, which has coefficients that converge
monotonically to zero from positive values,

\begin{displaymath}
  S_{z}
  =
  \sum_{k=1}^{\infty}
  \frac{1}{\e{\sqrt{k}}}\,
  z^{k}.
\end{displaymath}

\noindent
If we apply the ratio test to it, we get

\noindent
\begin{eqnarray*}
  R
  & = &
  \frac{\e{\sqrt{k}}|z|^{k+1}}{\e{\sqrt{k+1}}|z|^{k}}
  \\
  & = &
  \e{\sqrt{k}-\sqrt{k+1}}\rho
  \;\;\;\Rightarrow
  \\
  \ln(R)
  & = &
  \ln(\rho)
  +
  \sqrt{k}-\sqrt{k+1},
\end{eqnarray*}

\noindent
where $|z|=\rho$. In the large-$k$ limit we have

\noindent
\begin{eqnarray*}
  \lim_{k\to\infty}
  \ln(R)
  & = &
  \ln(\rho)
  +
  \lim_{k\to\infty}
  \left(
    \sqrt{k}-\sqrt{k+1}
  \right)
  \\
  & = &
  \ln(\rho)
  +
  \lim_{k\to\infty}
  \left[
    \sqrt{k}
    \left(
      1-\sqrt{1+\frac{1}{k}}
    \right)
  \right]
  \\
  & = &
  \ln(\rho)
  +
  \lim_{k\to\infty}
  \left[
    \sqrt{k}
    \left(
      1-1-\frac{1}{2k}
    \right)
  \right]
  \\
  & = &
  \ln(\rho)
  -
  \lim_{k\to\infty}
  \frac{1}{2\sqrt{k}}
  \\
  & = &
  \ln(\rho).
\end{eqnarray*}

\noindent
It follows that in the limit we have $R=\rho$ and therefore the conditions
of the test are satisfied if and only if $\rho<1$. This establishes the
open unit disk as the maximum disk of convergence of the series. Within
this disk the series converges to an inner analytic function $w(z)$, and
we may write

\begin{displaymath}
  w(z)
  =
  \sum_{k=1}^{\infty}
  \frac{1}{\e{\sqrt{k}}}\,
  z^{k}.
\end{displaymath}

\noindent
Since the maximum disk of convergence of the series is the open unit disk,
this function must have at least one singularity on the unit circle. Note
that since the series is monotonic with step $1$, we already know that
$w(z)$ has a single dominant singularity on that circle, located at $z=1$.
Consider now the logarithmic derivatives of this function. Using the
series we have within the open unit disk, for the $n^{\rm th}$ logarithmic
derivative of $w(z)$,

\begin{displaymath}
  w^{n\!\ldot}(z)
  =
  \sum_{k=1}^{\infty}
  \frac{k^{n}}{\e{\sqrt{k}}}\,
  z^{k}.
\end{displaymath}

\noindent
This notation includes the original series as the case $n=0$. All these
series converge on the open unit disk, of course. Let us now consider the
corresponding series of absolute values, for $\rho=1$,

\begin{displaymath}
  \overline{S}^{\,n\!\ldot}_{v}
  =
  \sum_{k=1}^{\infty}
  \frac{k^{n}}{\e{\sqrt{k}}}.
\end{displaymath}

\noindent
The terms of this sum can be bounded, for $k>k_{m}$ and some minimum value
$k_{m}$ of $k$, by the function $1/k^{2}$, since we have that it is always
possible to find a value $\xi_{m}$ of a variable $\xi$ such that for
$\xi>\xi_{m}$ we have

\begin{displaymath}
  \e{-\xi}
  <
  \frac{1}{\xi^{2n+4}},
\end{displaymath}

\noindent
since the exponential goes to zero faster than any inverse power. Making
$\xi^{2}=k$ and therefore $\xi=\sqrt{k}$ we have

\noindent
\begin{eqnarray*}
  \e{-\sqrt{k}}
  & < &
  \frac{1}{k^{n+2}}
  \;\;\;\Rightarrow
  \\
  \frac{k^{n}}{\e{\sqrt{k}}}
  & < &
  \frac{1}{k^{2}},
\end{eqnarray*}

\noindent
thus proving the assertion. This implies that all these series are
absolutely and uniformly convergent over the whole unit circle, to
continuous functions. Therefore, all these series must have a soft
singularity on the unit circle, at $z=1$. This is one example in which we
may differentiate as many times as we will, without the singularity ever
becoming hard. Therefore, that singularity is necessarily an infinitely
soft one.

Note that in this case the DP Fourier series on the unit circle converge
to $C^{\infty}$ functions, although there are singularities on that
circle. Although these real functions are $C^{\infty}$, in the real sense
of this concept, the complex function $w(z)$ cannot be $C^{\infty}$ on the
unit circle, in the complex sense of the concept. One may ask how can a
restriction of the complex function $w(z)$, which is $C^{\infty}$ in the
whole interior of the unit disk, be $C^{\infty}$ at the boundary of the
disk while $w(z)$ itself is not.

The answer is that the $C^{\infty}$ condition in the real sense is a
weaker condition than the $C^{\infty}$ condition in the complex sense.
While in the case of the real functions on the unit circle only the
derivatives with respect to $\theta$ must exist, in the case of the
complex function the derivatives in the perpendicular direction, that is
those with respect to $\rho$, must also exist, and in fact must give the
same values as the derivatives in the direction of $\theta$.

Unlike real functions over one-dimensional domains, which can be folded
around at will, complex analytic functions over two-dimensional domains
are rigid objects. If one restricts such a function to a one-dimensional
domain and then folds that domain around, the resulting real function over
it may no longer be the restriction of a complex function to the new
one-dimensional domain resulting from the folding process.

\section{Appendix: Examples of Center Series}\label{APPcenter}

In this appendix we will give a few simple illustrative examples of the
construction of center series. In many cases the construction of a center
series constitutes a practical way to determine the corresponding real
function, and to thus take advantage of the very existence of the inner
analytic function, for example when it is not possible to exhibit that
function in closed form, in order to explicitly take its limit to the unit
circle.

The process of construction consists of three parts, starting with the
determination of the power series $S_{z}$ from the original DP Fourier
series, which is simple and can always be done without any difficulty. The
second step is the construction of the complex center series $C_{z}$,
which is operationally fairly simple but depends on the knowledge of the
complete set of dominant singularities over the unit circle of the inner
analytic function $w(z)$ that the series $S_{z}$ converges to.

The last step is the recovery from $S_{z}$ written in terms of $C_{z}$ of
the real and imaginary parts of $w(z)$, in order to obtain the center
series versions of the original DP Fourier series and of its FC series.
This is straightforward but can become, in some cases, a rather long
algebraic process. In each example we will develop explicitly all these
steps, with a reasonable amount of detail.

Some of the examples that follow are the same that were worked out by
another method in the appendices of the already mentioned previous
paper~\cite{FTotCPI}. They are presented in the same order as in that
paper. It is understood that the final forms obtained for the functions
$f_{\rm s}(\theta)$ and $f_{\rm c}(\theta)$ in terms of the center series
are valid only away from the special points.

\subsection{A Regular Sine Series with All {\boldmath $k$}}

Consider the Fourier series of the one-cycle unit-amplitude sawtooth wave,
which is just the linear function $\theta/\pi$ between $-\pi$ and $\pi$.
As is well known it is given by the sine series

\begin{displaymath}
  S_{\rm s}
  =
  -\,
  \frac{2}{\pi}
  \sum_{k=1}^{\infty}
  \frac{(-1)^{k}}{k}\,
  \sin(k\theta).
\end{displaymath}

\noindent
The corresponding FC series is then

\begin{displaymath}
  \bar{S}_{\rm s}
  =
  -\,
  \frac{2}{\pi}
  \sum_{k=1}^{\infty}
  \frac{(-1)^{k}}{k}\,
  \cos(k\theta),
\end{displaymath}

\noindent
the complex $S_{v}$ series is given by

\begin{displaymath}
  S_{v}
  =
  -\,
  \frac{2}{\pi}
  \sum_{k=1}^{\infty}
  \frac{(-1)^{k}}{k}\,
  v^{k},
\end{displaymath}

\noindent
and the complex power series $S_{z}$ is given by

\begin{displaymath}
  S_{z}
  =
  -\,
  \frac{2}{\pi}
  \sum_{k=1}^{\infty}
  \frac{(-1)^{k}}{k}\,
  z^{k}.
\end{displaymath}

\noindent
The ratio test tells us that the disk of convergence of $S_{z}$ is the
unit disk. If we consider the inner analytic function $w(z)$ within this
disk we observe that $w(0)=0$, as expected. We have for this function

\begin{displaymath}
  w(z)
  =
  -\,
  \frac{2}{\pi}
  \sum_{k=1}^{\infty}
  \frac{(-1)^{k}}{k}\,
  z^{k}.
\end{displaymath}

\noindent
Being given by a monotonic series of step $1$ modified by the factor of
$(-1)^{k}$, this function has a single dominant singularity at $z=-1$,
where it diverges to infinity, as one can easily verify

\noindent
\begin{eqnarray*}
  w(-1)
  & = &
  -\,
  \frac{2}{\pi}
  \sum_{k=1}^{\infty}
  \frac{1}{k}
  \\
  & \to &
  -\infty.
\end{eqnarray*}

\noindent
We must therefore use a single factor of $(z+1)$ in the construction of
the center series,

\noindent
\begin{eqnarray*}
  C_{z}
  & = &
  -\,
  \frac{2}{\pi}\,
  (z+1)
  \sum_{k=1}^{\infty}
  \frac{(-1)^{k}}{k}\,
  z^{k}
  \\
  & = &
  -\,
  \frac{2}{\pi}
  \left[
    \,
    \sum_{k=1}^{\infty}
    \frac{(-1)^{k}}{k}\,
    z^{k+1}
    +
    \sum_{k=1}^{\infty}
    \frac{(-1)^{k}}{k}\,
    z^{k}
  \right]
  \\
  & = &
  -\,
  \frac{2}{\pi}
  \left[
    \,
    \sum_{k=2}^{\infty}
    \frac{(-1)^{k-1}}{k-1}\,
    z^{k}
    -
    z
    +
    \sum_{k=2}^{\infty}
    \frac{(-1)^{k}}{k}\,
    z^{k}
  \right]
  \\
  & = &
  -\,
  \frac{2}{\pi}
  \left[
    -z
    -
    \sum_{k=2}^{\infty}
    (-1)^{k}
    \left(
      \frac{1}{k-1}
      -
      \frac{1}{k}
    \right)
    z^{k}
  \right]
  \\
  & = &
  \frac{2}{\pi}
  \left[
    z
    +
    \sum_{k=2}^{\infty}
    \frac{(-1)^{k}}{(k-1)k}\,
    z^{k}
  \right]
  \\
  & = &
  \frac{2}{\pi}
  \left[
    z
    +
    \sum_{k=1}^{\infty}
    \frac{(-1)^{k+1}}{k(k+1)}\,
    z^{k+1}
  \right]
  \\
  & = &
  \frac{2}{\pi}\,
  z
  \left[
    1
    -
    \sum_{k=1}^{\infty}
    \frac{(-1)^{k}}{k(k+1)}\,
    z^{k}
  \right].
\end{eqnarray*}

\noindent
Unlike the original series, with coefficients that behave as $1/k$, this
series has coefficients that go to zero as $1/k^{2}$ when $k\to\infty$,
and therefore is absolutely and uniformly convergent to a continuous
function. This shows, in particular, that our evaluation of the set of
dominant singularities of $w(z)$ was in fact correct. We have therefore
for $S_{z}$ the representation

\begin{displaymath}
  S_{z}
  =
  \frac{2}{\pi}\,
  \frac{z}{z+1}
  \left[
    1
    -
    \sum_{k=1}^{\infty}
    \frac{(-1)^{k}}{k(k+1)}\,
    z^{k}
  \right],
\end{displaymath}

\noindent
with the singularity factored out and where the series involved is
absolutely and uniformly convergent, and therefore converges much faster
than the original one.

We may now take the real and imaginary parts of the $S_{v}$ series in
order to obtain faster-converging representation of the original DP
Fourier series and its FC series. We have on the unit circle

\noindent
\begin{eqnarray*}
  \frac{z}{z+1}
  & = &
  \frac{z(z^{*}+1)}{(z+1)(z^{*}+1)}
  \\
  & = &
  \frac{1+z}{2+z+z^{*}}
  \\
  & = &
  \frac{1+\cos(\theta)+\ii\sin(\theta)}{2+2\cos(\theta)}
  \\
  & = &
  \frac{1}{2}
  +
  \frac{\ii}{2}\,
  \frac{\sin(\theta)}{1+\cos(\theta)}.
\end{eqnarray*}

\noindent
If we write this in terms of $\theta/2$ we get

\noindent
\begin{eqnarray*}
  \frac{z}{z+1}
  & = &
  \frac{1}{2}
  +
  \frac{\ii}{2}\,
  \frac{2\sin(\theta/2)\cos(\theta/2)}{2\cos^{2}(\theta/2)}
  \\
  & = &
  \frac{1}{2}
  +
  \frac{\ii}{2}\,
  \frac{\sin(\theta/2)}{\cos(\theta/2)},
\end{eqnarray*}

\noindent
and therefore

\noindent
\begin{eqnarray*}
  S_{v}
  & = &
  \frac{1}{\pi}
  \left[
    1
    +
    \ii\,
    \frac{\sin(\theta/2)}{\cos(\theta/2)}
  \right]
  \left[
    1
    -
    \sum_{k=1}^{\infty}
    \frac{(-1)^{k}}{k(k+1)}\,
    \cos(k\theta)
    -
    \ii
    \sum_{k=1}^{\infty}
    \frac{(-1)^{k}}{k(k+1)}\,
    \sin(k\theta)
  \right]
  \\
  & = &
  \frac{1}{\pi}
  \left\{
    \left[
      1
      -
      \sum_{k=1}^{\infty}
      \frac{(-1)^{k}}{k(k+1)}\,
      \cos(k\theta)
    \right]
    +
    \frac{\sin(\theta/2)}{\cos(\theta/2)}
    \left[
      \,
      \sum_{k=1}^{\infty}
      \frac{(-1)^{k}}{k(k+1)}\,
      \sin(k\theta)
    \right]
  \right\}
  +
  \\
  &   &
  +
  \ii\,
  \frac{1}{\pi}
  \left\{
    \frac{\sin(\theta/2)}{\cos(\theta/2)}
    \left[
      1
      -
      \sum_{k=1}^{\infty}
      \frac{(-1)^{k}}{k(k+1)}\,
      \cos(k\theta)
    \right]
    -
    \left[
      \,
      \sum_{k=1}^{\infty}
      \frac{(-1)^{k}}{k(k+1)}\,
      \sin(k\theta)
    \right]
  \right\}
  \\
  & = &
  \frac{1}{\pi}
  \left\{
    1
    -
    \frac{1}{\cos(\theta/2)}
    \sum_{k=1}^{\infty}
    \frac{(-1)^{k}}{k(k+1)}
    \left[
      \cos\!\left(\frac{\theta}{2}\right)
      \cos(k\theta)
      -
      \sin\!\left(\frac{\theta}{2}\right)
      \sin(k\theta)
    \right]
  \right\}
  +
  \\
  &   &
  +
  \ii\,
  \frac{1}{\pi}
  \left\{
    \frac{\sin(\theta/2)}{\cos(\theta/2)}
    -
    \frac{1}{\cos(\theta/2)}
    \sum_{k=1}^{\infty}
    \frac{(-1)^{k}}{k(k+1)}
    \left[
      \sin\!\left(\frac{\theta}{2}\right)
      \cos(k\theta)
      +
      \cos\!\left(\frac{\theta}{2}\right)
      \sin(k\theta)
    \right]
  \right\}
  \\
  & = &
  \frac{1}{\pi\cos(\theta/2)}
  \left[
    \cos(\theta/2)
    -
    \sum_{k=1}^{\infty}
    \frac{(-1)^{k}}{k(k+1)}\,
    \cos\!\left(\frac{2k+1}{2}\,\theta\right)
  \right]
  +
  \\
  &   &
  +
  \ii\,
  \frac{1}{\pi\cos(\theta/2)}
  \left[
    \sin(\theta/2)
    -
    \sum_{k=1}^{\infty}
    \frac{(-1)^{k}}{k(k+1)}\,
    \sin\!\left(\frac{2k+1}{2}\,\theta\right)
  \right].
\end{eqnarray*}

\noindent
The original DP function is given by the imaginary part,

\begin{displaymath}
  f_{\rm s}(\theta)
  =
  \frac{1}{\pi\cos(\theta/2)}
  \left[
    \sin(\theta/2)
    -
    \sum_{k=1}^{\infty}
    \frac{(-1)^{k}}{k(k+1)}\,
    \sin\!\left(\frac{2k+1}{2}\,\theta\right)
  \right],
\end{displaymath}

\noindent
and the corresponding FC function $f_{\rm c}(\theta)=\bar{f}_{\rm
  s}(\theta)$ is given by the real part,

\begin{displaymath}
  f_{\rm c}(\theta)
  =
  \frac{1}{\pi\cos(\theta/2)}
  \left[
    \cos(\theta/2)
    -
    \sum_{k=1}^{\infty}
    \frac{(-1)^{k}}{k(k+1)}\,
    \cos\!\left(\frac{2k+1}{2}\,\theta\right)
  \right].
\end{displaymath}

\noindent
Both of these series are absolutely and uniformly convergent.

\subsection{A Regular Sine Series with Odd {\boldmath
    $k$}}\label{APPregsqwave}

Consider the Fourier series of the standard unit-amplitude square wave. As
is well known it is given by the sine series

\begin{displaymath}
  S_{\rm s}
  =
  \frac{4}{\pi}
  \sum_{j=0}^{\infty}
  \frac{1}{2j+1}\,
  \sin[(2j+1)\theta].
\end{displaymath}

\noindent
The corresponding FC series is then

\begin{displaymath}
  \bar{S}_{\rm s}
  =
  \frac{4}{\pi}
  \sum_{j=0}^{\infty}
  \frac{1}{2j+1}\,
  \cos[(2j+1)\theta],
\end{displaymath}

\noindent
the complex $S_{v}$ series is given by

\begin{displaymath}
  S_{v}
  =
  \frac{4}{\pi}
  \sum_{j=0}^{\infty}
  \frac{1}{2j+1}\,
  v^{2j+1},
\end{displaymath}

\noindent
and the complex power series $S_{z}$ is given by

\begin{displaymath}
  S_{z}
  =
  \frac{4}{\pi}
  \sum_{j=0}^{\infty}
  \frac{1}{2j+1}\,
  z^{2j+1}.
\end{displaymath}

\noindent
The ratio test tells us that the disk of convergence of $S_{z}$ is the
unit disk. If we consider the inner analytic function $w(z)$ within this
disk we observe that $w(0)=0$, as expected. We have for this function

\begin{displaymath}
  w(z)
  =
  \frac{4}{\pi}
  \sum_{j=0}^{\infty}
  \frac{1}{2j+1}\,
  z^{2j+1}.
\end{displaymath}

\noindent
Being given by a monotonic series of step $2$ this function has two
dominant singularities, located at $z=1$ and at $z=-1$, where it diverges
to infinity, as one can easily verify,

\noindent
\begin{eqnarray*}
  w(1)
  & = &
  \frac{4}{\pi}
  \sum_{j=0}^{\infty}
  \frac{1}{2j+1}
  \\
  & \to &
  \infty,
  \\
  w(-1)
  & = &
  -\,
  \frac{4}{\pi}
  \sum_{j=0}^{\infty}
  \frac{1}{2j+1}
  \\
  & \to &
  -\infty.
\end{eqnarray*}

\noindent
We must therefore use the two factors $(z-1)(z+1)=z^{2}-1$ in the
construction of the center series,

\noindent
\begin{eqnarray*}
  C_{z}
  & = &
  \frac{4}{\pi}\,
  \left(z^{2}-1\right)
  \sum_{j=0}^{\infty}
  \frac{1}{2j+1}\,
  z^{2j+1}
  \\
  & = &
  \frac{4}{\pi}
  \left(
    \sum_{j=0}^{\infty}
    \frac{1}{2j+1}\,
    z^{2j+3}
    -
    \sum_{j=0}^{\infty}
    \frac{1}{2j+1}\,
    z^{2j+1}
  \right)
  \\
  & = &
  \frac{4}{\pi}
  \left(
    \sum_{j=1}^{\infty}
    \frac{1}{2j-1}\,
    z^{2j+1}
    -
    z
    -
    \sum_{j=1}^{\infty}
    \frac{1}{2j+1}\,
    z^{2j+1}
  \right)
  \\
  & = &
  \frac{4}{\pi}
  \left[
    -z
    +
    \sum_{j=1}^{\infty}
    \left(
      \frac{1}{2j-1}
      -
      \frac{1}{2j+1}
    \right)
    z^{2j+1}
  \right]
  \\
  & = &
  \frac{4}{\pi}
  \left(
    -z
    +
    \sum_{j=1}^{\infty}
    \frac{2}{4j^{2}-1}\,
    z^{2j+1}
  \right)
  \\
  & = &
  \frac{4}{\pi}\,
  z
  \left(
    -1
    +
    \sum_{j=1}^{\infty}
    \frac{2}{4j^{2}-1}\,
    z^{2j}
  \right).
\end{eqnarray*}

\noindent
Unlike the original series, with coefficients that behave as $1/k$ (with
$k=2j+1$), this series has coefficients that go to zero as $1/k^{2}$ when
$k\to\infty$, and therefore is absolutely and uniformly convergent to a
continuous function. This shows, in particular, that our evaluation of the
set of dominant singularities of $w(z)$ was in fact correct. We have
therefore for $S_{z}$ the representation

\begin{displaymath}
  S_{z}
  =
  \frac{4}{\pi}\,
  \frac{z}{z^{2}-1}
  \left(
    -1
    +
    \sum_{j=1}^{\infty}
    \frac{2}{4j^{2}-1}\,
    z^{2j}
  \right),
\end{displaymath}

\noindent
with the singularities factored out and where the series involved is
absolutely and uniformly convergent, and therefore converges much faster
than the original one.

We may now take the real and imaginary parts of the $S_{v}$ series in
order to obtain faster-converging representation of the original DP
Fourier series and its FC series. We have on the unit circle

\noindent
\begin{eqnarray*}
  \frac{z}{z^{2}-1}
  & = &
  \frac
  {z\left[(z^{*})^{2}-1\right]}
  {\left(z^{2}-1\right)\left[(z^{*})^{2}-1\right]}
  \\
  & = &
  \frac{z^{*}-z}{2-z^{2}-(z^{*})^{2}}
  \\
  & = &
  \frac{-2\ii\sin(\theta)}{2-2\cos(2\theta)}
  \\
  & = &
  \frac{-\ii\sin(\theta)}{1-\cos^{2}(\theta)+\sin^{2}(\theta)}
  \\
  & = &
  \frac{-\ii\sin(\theta)}{2\sin^{2}(\theta)}
  \\
  & = &
  \frac{-\ii}{2\sin(\theta)},
\end{eqnarray*}

\noindent
and therefore

\noindent
\begin{eqnarray*}
  S_{v}
  & = &
  \frac{4}{\pi}\,
  \frac{-\ii}{2\sin(\theta)}
  \left[
    -1
    +
    \sum_{j=1}^{\infty}
    \frac{2}{4j^{2}-1}\,
    \cos(2j\theta)
    +
    \ii
    \sum_{j=1}^{\infty}
    \frac{2}{4j^{2}-1}\,
    \sin(2j\theta)
  \right]
  \\
  & = &
  \frac{2}{\pi\sin(\theta)}
  \left[
    \sum_{j=1}^{\infty}
    \frac{2}{4j^{2}-1}\,
    \sin(2j\theta)
  \right]    
  +
  \ii\,
  \frac{2}{\pi\sin(\theta)}
  \left[
    1
    -
    \sum_{j=1}^{\infty}
    \frac{2}{4j^{2}-1}\,
    \cos(2j\theta)
  \right].
\end{eqnarray*}

\noindent
The original DP function is given by the imaginary part,

\begin{displaymath}
  f_{\rm s}(\theta)
  =
  \frac{2}{\pi\sin(\theta)}
  \left[
    1
    -
    \sum_{j=1}^{\infty}
    \frac{2}{4j^{2}-1}\,
    \cos(2j\theta)
  \right],
\end{displaymath}

\noindent
and the corresponding FC function $f_{\rm c}(\theta)=\bar{f}_{\rm
  s}(\theta)$ is given by the real part,

\begin{displaymath}
  f_{\rm c}(\theta)
  =
  \frac{2}{\pi\sin(\theta)}
  \left[
    \sum_{j=1}^{\infty}
    \frac{2}{4j^{2}-1}\,
    \sin(2j\theta)
  \right].
\end{displaymath}

\noindent
Both of these series are absolutely and uniformly convergent.

\subsection{A Regular Sine Series with Even {\boldmath $k$}}

Consider the Fourier series of the two-cycle unit-amplitude sawtooth wave.
As is well known it is given by the sine series

\begin{displaymath}
  S_{\rm s}
  =
  -\,
  \frac{4}{\pi}
  \sum_{j=1}^{\infty}
  \frac{1}{2j}\,
  \sin[(2j)\theta].
\end{displaymath}

\noindent
The corresponding FC series is then

\begin{displaymath}
  \bar{S}_{\rm s}
  =
  -\,
  \frac{4}{\pi}
  \sum_{j=1}^{\infty}
  \frac{1}{2j}\,
  \cos[(2j)\theta],
\end{displaymath}

\noindent
the complex $S_{v}$ series is given by

\begin{displaymath}
  S_{v}
  =
  -\,
  \frac{4}{\pi}
  \sum_{j=1}^{\infty}
  \frac{1}{2j}\,
  v^{2j},
\end{displaymath}

\noindent
and the complex power series $S_{z}$ is given by

\begin{displaymath}
  S_{z}
  =
  -\,
  \frac{4}{\pi}
  \sum_{j=1}^{\infty}
  \frac{1}{2j}\,
  z^{2j}.
\end{displaymath}

\noindent
The ratio test tells us that the disk of convergence of $S_{z}$ is the
unit disk. If we consider the inner analytic function $w(z)$ within this
disk we observe that $w(0)=0$, as expected. We have for this function

\begin{displaymath}
  w(z)
  =
  -\,
  \frac{4}{\pi}
  \sum_{j=1}^{\infty}
  \frac{1}{2j}\,
  z^{2j}.
\end{displaymath}

\noindent
Being given by a monotonic series of step $2$ this function has two
dominant singularities, located at $z=1$ and at $z=-1$, where it diverges
to infinity, as one can easily verify,

\noindent
\begin{eqnarray*}
  w(1)
  & = &
  -\,
  \frac{4}{\pi}
  \sum_{j=1}^{\infty}
  \frac{1}{2j}\,
  \\
  & \to &
  -\infty,
  \\
  w(-1)
  & = &
  -\,
  \frac{4}{\pi}
  \sum_{j=1}^{\infty}
  \frac{1}{2j}\,
  \\
  & \to &
  -\infty.
\end{eqnarray*}

\noindent
We must therefore use the two factors $(z-1)(z+1)=z^{2}-1$ in the
construction of the center series,

\noindent
\begin{eqnarray*}
  C_{z}
  & = &
  -\,
  \frac{4}{\pi}\,
  \left(z^{2}-1\right)
  \sum_{j=1}^{\infty}
  \frac{1}{2j}\,
  z^{2j}
  \\
  & = &
  -\,
  \frac{4}{\pi}
  \left(
    \sum_{j=1}^{\infty}
    \frac{1}{2j}\,
    z^{2j+2}
    -
    \sum_{j=1}^{\infty}
    \frac{1}{2j}\,
    z^{2j}
  \right)
  \\
  & = &
  -\,
  \frac{4}{\pi}
  \left(
    \sum_{j=2}^{\infty}
    \frac{1}{2j-2}\,
    z^{2j}
    -
    \frac{z^{2}}{2}
    -
    \sum_{j=2}^{\infty}
    \frac{1}{2j}\,
    z^{2j}
  \right)
  \\
  & = &
  -\,
  \frac{2}{\pi}
  \left[
    -z^{2}
    +
    \sum_{j=2}^{\infty}
    \left(
      \frac{1}{j-1}
      -
      \frac{1}{j}
    \right)
    z^{2j}
  \right]
  \\
  & = &
  \frac{2}{\pi}
  \left[
    z^{2}
    -
    z^{2}
    \sum_{j=2}^{\infty}
    \frac{1}{(j-1)j}\,
    z^{2j-2}
  \right]
  \\
  & = &
  \frac{2}{\pi}\,
  z^{2}
  \left[
    1
    -
    \sum_{j=1}^{\infty}
    \frac{1}{j(j+1)}\,
    z^{2j}
  \right].
\end{eqnarray*}

\noindent
Unlike the original series, with coefficients that behave as $1/k$ (with
$k=2j$), this series has coefficients that go to zero as $1/k^{2}$ when
$k\to\infty$, and therefore is absolutely and uniformly convergent to a
continuous function. This shows, in particular, that our evaluation of the
set of dominant singularities of $w(z)$ was in fact correct. We have
therefore for $S_{z}$ the representation

\begin{displaymath}
  S_{z}
  =
  \frac{2}{\pi}\,
  \frac{z^{2}}{z^{2}-1}
  \left[
    1
    -
    \sum_{j=1}^{\infty}
    \frac{1}{j(j+1)}\,
    z^{2j}
  \right],
\end{displaymath}

\noindent
with the singularities factored out and where the series involved is
absolutely and uniformly convergent, and therefore converges much faster
than the original one.

We may now take the real and imaginary parts of the $S_{v}$ series in
order to obtain faster-converging representation of the original DP
Fourier series and its FC series. We have on the unit circle

\noindent
\begin{eqnarray*}
  \frac{z^{2}}{z^{2}-1}
  & = &
  \frac
  {z^{2}\left[(z^{*})^{2}-1\right]}
  {\left(z^{2}-1\right)\left[(z^{*})^{2}-1\right]}
  \\
  & = &
  \frac{1-z^{2}}{2-z^{2}-(z^{*})^{2}}
  \\
  & = &
  \frac{1-\cos(2\theta)-\ii\sin(2\theta)}{2-2\cos(2\theta)}
  \\
  & = &
  \frac{1}{2}
  -
  \frac{\ii}{2}\,
  \frac{\sin(2\theta)}{1-\cos(2\theta)}
  \\
  & = &
  \frac{1}{2}
  -
  \frac{\ii}{2}\,
  \frac{2\sin(\theta)\cos(\theta)}{2\sin^{2}(\theta)}
  \\
  & = &
  \frac{1}{2}
  -
  \frac{\ii}{2}\,
  \frac{\cos(\theta)}{\sin(\theta)},
\end{eqnarray*}

\noindent
and therefore

\noindent
\begin{eqnarray*}
  S_{v}
  & = &
  \frac{1}{\pi}
  \left[
    1
    -
    \ii\,
    \frac{\cos(\theta)}{\sin(\theta)}
  \right]
  \left[
    1
    -
    \sum_{j=1}^{\infty}
    \frac{1}{j(j+1)}\,
    \cos(2j\theta)
    -
    \ii
    \sum_{j=1}^{\infty}
    \frac{1}{j(j+1)}\,
    \sin(2j\theta)
  \right]
  \\
  & = &
  \frac{1}{\pi}
  \left\{
    \left[
      1
      -
      \sum_{j=1}^{\infty}
      \frac{1}{j(j+1)}\,
      \cos(2j\theta)
    \right]
    -
    \frac{\cos(\theta)}{\sin(\theta)}
    \left[
      \sum_{j=1}^{\infty}
      \frac{1}{j(j+1)}\,
      \sin(2j\theta)
    \right]
  \right\}
  +
  \\
  &   &
  +
  \ii\,
  \frac{1}{\pi}
  \left\{
    -\,
    \frac{\cos(\theta)}{\sin(\theta)}
    \left[
      1
      -
      \sum_{j=1}^{\infty}
      \frac{1}{j(j+1)}\,
      \cos(2j\theta)
    \right]
    -
    \left[
      \sum_{j=1}^{\infty}
      \frac{1}{j(j+1)}\,
      \sin(2j\theta)
    \right]
  \right\}
  \\
  & = &
  \frac{1}{\pi}
  \left\{
    1
    -
    \frac{1}{\sin(\theta)}
    \sum_{j=1}^{\infty}
    \frac{1}{j(j+1)}
    \left[
      \sin(\theta)
      \cos(2j\theta)
      +
      \cos(\theta)
      \sin(2j\theta)
    \right]
  \right\}
  +
  \\
  &   &
  +
  \ii\,
  \frac{1}{\pi}
  \left\{
    -\,
    \frac{\cos(\theta)}{\sin(\theta)}
    +
    \frac{1}{\sin(\theta)}
    \sum_{j=1}^{\infty}
    \frac{1}{j(j+1)}
    \left[
      \cos(\theta)
      \cos(2j\theta)
      -
      \sin(\theta)
      \sin(2j\theta)
    \right]
  \right\}
  \\
  & = &
  \frac{1}{\pi\sin(\theta)}
  \left\{
    \sin(\theta)
    -
    \sum_{j=1}^{\infty}
    \frac{1}{j(j+1)}\,
    \sin[(2j+1)\theta]
  \right\}
  +
  \\
  &   &
  +
  \ii\,
  \frac{1}{\pi\sin(\theta)}
  \left\{
    -
    \cos(\theta)
    +
    \sum_{j=1}^{\infty}
    \frac{1}{j(j+1)}\,
    \cos[(2j+1)\theta]
  \right\}.
\end{eqnarray*}

\noindent
The original DP function is given by the imaginary part,

\begin{displaymath}
  f_{\rm s}(\theta)
  =
  -\,
  \frac{1}{\pi\sin(\theta)}
  \left\{
    \cos(\theta)
    -
    \sum_{j=1}^{\infty}
    \frac{1}{j(j+1)}\,
    \cos[(2j+1)\theta]
  \right\},
\end{displaymath}

\noindent
and the corresponding FC function $f_{\rm c}(\theta)=\bar{f}_{\rm
  s}(\theta)$ is given by the real part,

\begin{displaymath}
  f_{\rm c}(\theta)
  =
  \frac{1}{\pi\sin(\theta)}
  \left\{
    \sin(\theta)
    -
    \sum_{j=1}^{\infty}
    \frac{1}{j(j+1)}\,
    \sin[(2j+1)\theta]
  \right\}.
\end{displaymath}

\noindent
Both of these series are absolutely and uniformly convergent.

\subsection{A Regular Cosine Series with Odd {\boldmath $k$}}

Consider the Fourier series of the unit-amplitude triangular wave. As is
well known it is given by the cosine series

\begin{displaymath}
  S_{\rm c}
  =
  -\,
  \frac{8}{\pi^{2}}
  \sum_{j=0}^{\infty}
  \frac{1}{(2j+1)^{2}}\,
  \cos[(2j+1)\theta].
\end{displaymath}

\noindent
The corresponding FC series is then

\begin{displaymath}
  \bar{S}_{\rm c}
  =
  -\,
  \frac{8}{\pi^{2}}
  \sum_{j=0}^{\infty}
  \frac{1}{(2j+1)^{2}}\,
  \sin[(2j+1)\theta].
\end{displaymath}

\noindent
Note that due to the factors of $1/k^{2}$ (with $k=2j+1$), these series
are already absolutely and uniformly convergent. But we will proceed with
the construction in any case. The complex $S_{v}$ series is given by

\begin{displaymath}
  S_{v}
  =
  -\,
  \frac{8}{\pi^{2}}
  \sum_{j=0}^{\infty}
  \frac{1}{(2j+1)^{2}}\,
  v^{2j+1},
\end{displaymath}

\noindent
and the complex power series $S_{z}$ is given by

\begin{displaymath}
  S_{z}
  =
  -\,
  \frac{8}{\pi^{2}}
  \sum_{j=0}^{\infty}
  \frac{1}{(2j+1)^{2}}\,
  z^{2j+1}.
\end{displaymath}

\noindent
The ratio test tells us that the disk of convergence of $S_{z}$ is the
unit disk. If we consider the inner analytic function $w(z)$ within this
disk we observe that $w(0)=0$, as expected. We have for this function

\begin{displaymath}
  w(z)
  =
  -\,
  \frac{8}{\pi^{2}}
  \sum_{j=0}^{\infty}
  \frac{1}{(2j+1)^{2}}\,
  z^{2j+1}.
\end{displaymath}

\noindent
Being given by a monotonic series of step $2$ this function has two
dominant singularities, located at $z=1$ and at $z=-1$, as one can easily
verify by taking its logarithmic derivative, which is proportional to the
inner analytic function of the square wave, that we examined before in
Subsection~\ref{APPregsqwave},

\begin{displaymath}
  z\,
  \frac{dw(z)}{dz}
  =
  -\,
  \frac{8}{\pi^{2}}
  \sum_{j=0}^{\infty}
  \frac{1}{2j+1}\,
  z^{2j+1}.
\end{displaymath}

\noindent
We must therefore use the two factors $(z-1)(z+1)=z^{2}-1$ in the
construction of the center series,

\noindent
\begin{eqnarray*}
  C_{z}
  & = &
  -\,
  \frac{8}{\pi^{2}}\,
  \left(z^{2}-1\right)
  \sum_{j=0}^{\infty}
  \frac{1}{(2j+1)^{2}}\,
  z^{2j+1}
  \\
  & = &
  -\,
  \frac{8}{\pi^{2}}
  \left[
    \sum_{j=0}^{\infty}
    \frac{1}{(2j+1)^{2}}\,
    z^{2j+3}
    -
    \sum_{j=0}^{\infty}
    \frac{1}{(2j+1)^{2}}\,
    z^{2j+1}
  \right]
  \\
  & = &
  -\,
  \frac{8}{\pi^{2}}
  \left[
    \sum_{j=1}^{\infty}
    \frac{1}{(2j-1)^{2}}\,
    z^{2j+1}
    -
    z
    -
    \sum_{j=1}^{\infty}
    \frac{1}{(2j+1)^{2}}\,
    z^{2j+1}
  \right]
  \\
  & = &
  -\,
  \frac{8}{\pi^{2}}
  \left\{
    -z
    +
    \sum_{j=1}^{\infty}
    \left[
      \frac{1}{(2j-1)^{2}}
      -
      \frac{1}{(2j+1)^{2}}
    \right]
    z^{2j+1}
  \right\}
  \\
  & = &
  \frac{8}{\pi^{2}}
  \left[
    z
    -
    z
    \sum_{j=1}^{\infty}
    \frac{8j}{\left(4j^{2}-1\right)^{2}}\,
    z^{2j}
  \right]
  \\
  & = &
  \frac{8}{\pi^{2}}\,
  z
  \left[
    1
    -
    \sum_{j=1}^{\infty}
    \frac{8j}{\left(4j^{2}-1\right)^{2}}\,
    z^{2j}
  \right].
\end{eqnarray*}

\noindent
Unlike the original series, with coefficients that behave as $1/k^{2}$
(with $k=2j+1$), this series has coefficients that go to zero as $1/k^{3}$
when $k\to\infty$, and therefore converges faster than the original one.
This shows, in particular, that our evaluation of the set of dominant
singularities of $w(z)$ was in fact correct. We have therefore for $S_{z}$
the representation

\begin{displaymath}
  S_{z}
  =
  \frac{8}{\pi^{2}}\,
  \frac{z}{z^{2}-1}
  \left[
    1
    -
    \sum_{j=1}^{\infty}
    \frac{8j}{\left(4j^{2}-1\right)^{2}}\,
    z^{2j}
  \right],
\end{displaymath}

\noindent
with the singularities factored out. Although both this series and the
original one are absolutely and uniformly convergent, this converges
faster, and may be differentiated once, still resulting in another series
which is also absolutely and uniformly convergent. Note that in this case,
as was discussed in the appendices of the previous paper~\cite{FTotCPI},
we are not able to write an explicit expression for $w(z)$ in terms of
elementary function, so that we cannot explicitly take its limit to the
unit circle. However, as one can see here we are still able to write a
series to represent it over the unit circle, which converges faster that
the original one. This gives us the possibility of calculating the
function to any required precision level, and to do so efficiently.

We may now take the real and imaginary parts of the $S_{v}$ series in
order to obtain faster-converging representation of the original DP
Fourier series and its FC series. We have on the unit circle, as we saw
before in Subsection~\ref{APPregsqwave},

\begin{displaymath}
  \frac{z}{z^{2}-1}
  =
  \frac{-\ii}{2\sin(\theta)},
\end{displaymath}

\noindent
and therefore

\noindent
\begin{eqnarray*}
  S_{v}
  & = &
  \frac{4}{\pi^{2}}\,
  \frac{-\ii}{\sin(\theta)}
  \left[
    1
    -
    \sum_{j=1}^{\infty}
    \frac{8j}{\left(4j^{2}-1\right)^{2}}\,
    \cos(2j\theta)
    -
    \ii
    \sum_{j=1}^{\infty}
    \frac{8j}{\left(4j^{2}-1\right)^{2}}\,
    \sin(2j\theta)
  \right]
  \\
  & = &
  \frac{4}{\pi^{2}\sin(\theta)}
  \left[
    -
    \sum_{j=1}^{\infty}
    \frac{8j}{\left(4j^{2}-1\right)^{2}}\,
    \sin(2j\theta)
  \right]
  +
  \\
  &   &
  +
  \ii\,
  \frac{4}{\pi^{2}\sin(\theta)}
  \left[
    -1
    +
    \sum_{j=1}^{\infty}
    \frac{8j}{\left(4j^{2}-1\right)^{2}}\,
    \cos(2j\theta)
  \right].
\end{eqnarray*}

\noindent
The original DP function is given by the real part,

\begin{displaymath}
  f_{\rm c}(\theta)
  =
  -\,
  \frac{4}{\pi^{2}\sin(\theta)}
  \left[
    \sum_{j=1}^{\infty}
    \frac{8j}{\left(4j^{2}-1\right)^{2}}\,
    \sin(2j\theta)
  \right],
\end{displaymath}

\noindent
and the corresponding FC function $f_{\rm s}(\theta)=\bar{f}_{\rm
  c}(\theta)$ is given by the imaginary part,

\begin{displaymath}
  f_{\rm s}(\theta)
  =
  -\,
  \frac{4}{\pi^{2}\sin(\theta)}
  \left[
    1
    -
    \sum_{j=1}^{\infty}
    \frac{8j}{\left(4j^{2}-1\right)^{2}}\,
    \cos(2j\theta)
  \right].
\end{displaymath}

\noindent
Both of these series converge faster than the original Fourier series.

\subsection{A Singular Cosine Series}

Consider the Fourier series of the Dirac delta ``function'' centered at
$\theta=\theta_{1}$, which we denote by $\delta(\theta-\theta_{1})$. We
may easily calculate its Fourier coefficients using the rules of
manipulation of $\delta(\theta-\theta_{1})$, thus obtaining
$\alpha_{k}=\cos(k\theta_{1})/\pi$ and $\beta_{k}=\sin(k\theta_{1})/\pi$
for all $k$. The series is therefore the complete Fourier series given by

\noindent
\begin{eqnarray*}
  S_{\rm c}
  & = &
  \frac{1}{2\pi}
  +
  \frac{1}{\pi}
  \sum_{k=1}^{\infty}
  \left[
    \cos(k\theta_{1})
    \cos(k\theta)
    +
    \sin(k\theta_{1})
    \sin(k\theta)
  \right]
  \\
  & = &
  \frac{1}{2\pi}
  +
  \frac{1}{\pi}
  \sum_{k=1}^{\infty}
  \cos(k\Delta\theta),
\end{eqnarray*}

\noindent
where $\Delta\theta=\theta-\theta_{1}$. Apart from the constant term this
is in fact a DP cosine series on this new variable. Clearly, this series
diverges at all points in the interval $[-\pi,\pi]$. Undaunted by this, we
proceed to construct the FC series, with respect to the new variable
$\Delta\theta$, which turns out to be

\begin{displaymath}
  \bar{S}_{\rm c}
  =
  \frac{1}{\pi}
  \sum_{k=1}^{\infty}
  \sin(k\Delta\theta),
\end{displaymath}

\noindent
a series that is also divergent, this time almost everywhere. If we define
$v=\exp(\ii\theta)$ and $v_{1}=\exp(\ii\theta_{1})$ the corresponding
complex series $S_{v}$ is then given by

\begin{displaymath}
  S_{v}
  =
  \frac{1}{2\pi}
  +
  \frac{1}{\pi}
  \sum_{k=1}^{\infty}
  \left(\frac{v}{v_{1}}\right)^{k},
\end{displaymath}

\noindent
where we included the $k=0$ term, and the corresponding complex power
series $S_{z}$ is given by

\begin{displaymath}
  S_{z}
  =
  \frac{1}{2\pi}
  +
  \frac{1}{\pi}
  \sum_{k=1}^{\infty}
  \left(\frac{z}{z_{1}}\right)^{k},
\end{displaymath}

\noindent
where $z=\rho v$ and $z_{1}=v_{1}$ is a point over the unit circle. The
ratio test tells us that the disk of convergence of $S_{z}$ is the unit
disk. This converges to a perfectly well-defined analytic function
strictly inside the open unit disk. If we eliminate the constant term we
get a series $S'_{z}$ which converges to an inner analytic function
rotated by the angle $\theta_{1}$,

\begin{displaymath}
  S'_{z}
  =
  \frac{1}{\pi}
  \sum_{k=1}^{\infty}
  \left(\frac{z}{z_{1}}\right)^{k}.
\end{displaymath}

\noindent
The dominant singularity is clearly at the point $z_{1}$, so we must use
the factor $(z-z_{1})$ in the construction of the corresponding center
series,

\noindent
\begin{eqnarray*}
  C'_{z}
  & = &
  (z-z_{1})S'_{z}
  \\
  & = &
  \frac{1}{\pi}\,
  z_{1}
  \left(\frac{z}{z_{1}}-1\right)
  \sum_{k=1}^{\infty}
  \left(\frac{z}{z_{1}}\right)^{k}
  \\
  & = &
  \frac{1}{\pi}\,
  z_{1}
  \left[
    \,
    \sum_{k=1}^{\infty}
    \left(\frac{z}{z_{1}}\right)^{k+1}
    -
    \sum_{k=1}^{\infty}
    \left(\frac{z}{z_{1}}\right)^{k}
  \right]
  \\
  & = &
  \frac{1}{\pi}\,
  z_{1}
  \left[
    \,
    \sum_{k=2}^{\infty}
    \left(\frac{z}{z_{1}}\right)^{k}
    -
    \left(\frac{z}{z_{1}}\right)
    -
    \sum_{k=2}^{\infty}
    \left(\frac{z}{z_{1}}\right)^{k}
  \right]
  \\
  & = &
  -\,
  \frac{1}{\pi}\,
  z.
\end{eqnarray*}

\noindent
So we see that we get a remarkably simple result, since the center series
can actually be added up exactly. We get therefore for the series $S'_{z}$

\begin{displaymath}
  S'_{z}
  =
  -\,
  \frac{1}{\pi}\,
  \frac{z}{z-z_{1}},
\end{displaymath}

\noindent
and for the series $S_{z}$

\begin{displaymath}
  S_{z}
  =
  \frac{1}{2\pi}
  -
  \frac{1}{\pi}\,
  \frac{z}{z-z_{1}}.
\end{displaymath}

\noindent
We may now take the real and imaginary parts of the $S_{v}$ series in
order to obtain faster-converging representation of the original DP
Fourier series and its FC series. The explanation of the reasons why this
is a representation of the Dirac delta ``function'' requires taking limits
to the unit circle carefully, and since they were given in the previous
paper~\cite{FTotCPI}, they will not be repeated here. We have, for $z$ on
the unit circle, so long as $\Delta\theta\neq 0$,

\noindent
\begin{eqnarray*}
  \frac{z}{z-z_{1}}
  & = &
  \frac{z(z^{*}-z_{1}^{*})}{(z-z_{1})(z^{*}-z_{1}^{*})}
  \\
  & = &
  \frac{1-zz_{1}^{*}}{2-zz_{1}^{*}-z^{*}z_{1}}
  \\
  & = &
  \frac
  {1-\cos(\Delta\theta)-\ii\sin(\Delta\theta)}
  {2-2\cos(\Delta\theta)}
  \\
  & = &
  \frac{1}{2}
  -
  \ii\,
  \frac{1}{2}\,
  \frac{\sin(\Delta\theta)}{1-\cos(\Delta\theta)},
  \\
  & = &
  \frac{1}{2}
  -
  \ii\,
  \frac{1}{2}\,
  \frac{1+\cos(\Delta\theta)}{\sin(\Delta\theta)},
\end{eqnarray*}

\noindent
and therefore

\noindent
\begin{eqnarray*}
  S_{v}
  & = &
  \frac{1}{2\pi}
  -
  \frac{1}{2\pi}
  +
  \ii\,
  \frac{1}{2\pi}\,
  \frac{1+\cos(\Delta\theta)}{\sin(\Delta\theta)}
  \\
  & = &
  \ii\,
  \frac{1}{2\pi}\,
  \frac{1+\cos(\Delta\theta)}{\sin(\Delta\theta)}.
\end{eqnarray*}

\noindent
The original DP ``function'' is given by the real part, and therefore we
get

\begin{displaymath}
  f_{\rm c}(\theta)
  =
  0,
\end{displaymath}

\noindent
which is the correct value for the Dirac delta ``function'' away from the
singular point at $\Delta\theta=0$, and the corresponding FC function
$f_{\rm s}(\theta)=\bar{f}_{\rm c}(\theta)$ is given by the imaginary
part,

\begin{displaymath}
  f_{\rm s}(\theta)
  =
  \frac{1}{2\pi}\,
  \frac{1+\cos(\Delta\theta)}{\sin(\Delta\theta)},
\end{displaymath}

\noindent
which is the same result we obtained in the previous paper~\cite{FTotCPI}.

\subsection{Another Regular Cosine Series with Odd {\boldmath $k$}}

Consider the Fourier series of the unit-amplitude square wave, shifted
along the $\theta$ axis to $\theta'$, with $\theta-\theta'= \pi/2$, so
that it becomes an even function. As is well known it is given by the
cosine series

\begin{displaymath}
  S_{\rm c}
  =
  \frac{4}{\pi}
  \sum_{j=0}^{\infty}
  \frac{(-1)^{j}}{2j+1}\,
  \cos[(2j+1)\theta],
\end{displaymath}

\noindent
where we have dropped the prime. The corresponding FC series is then

\begin{displaymath}
  \bar{S}_{\rm c}
  =
  \frac{4}{\pi}
  \sum_{j=0}^{\infty}
  \frac{(-1)^{j}}{2j+1}\,
  \sin[(2j+1)\theta],
\end{displaymath}

\noindent
the complex $S_{v}$ series is given by

\begin{displaymath}
  S_{v}
  =
  \frac{4}{\pi}
  \sum_{j=0}^{\infty}
  \frac{(-1)^{j}}{2j+1}\,
  v^{2j+1},
\end{displaymath}

\noindent
and the complex power series $S_{z}$ is given by

\begin{displaymath}
  S_{z}
  =
  \frac{4}{\pi}
  \sum_{j=0}^{\infty}
  \frac{(-1)^{j}}{2j+1}\,
  z^{2j+1}.
\end{displaymath}

\noindent
The ratio test tells us that the disk of convergence of $S_{z}$ is the
unit disk. If we consider the inner analytic function $w(z)$ within this
disk we observe that $w(0)=0$, as expected. We have for this function

\begin{displaymath}
  w(z)
  =
  \frac{4}{\pi}
  \sum_{j=0}^{\infty}
  \frac{(-1)^{j}}{2j+1}\,
  z^{2j+1}.
\end{displaymath}

\noindent
Being given by a monotonic series of step $2$ modified by the factor of
$(-1)^{j}$, this function has two dominant singularities, located at
$z=\ii$ and at $z=-\ii$, where it diverges to infinity, as one can easily
verify,

\noindent
\begin{eqnarray*}
  w(\ii)
  & = &
  \frac{4\ii}{\pi}
  \sum_{j=0}^{\infty}
  \frac{1}{2j+1}
  \\
  & \to &
  \ii\infty,
  \\
  w(-\ii)
  & = &
  -\,
  \frac{4\ii}{\pi}
  \sum_{j=0}^{\infty}
  \frac{1}{2j+1}
  \\
  & \to &
  -\ii\infty.
\end{eqnarray*}

\noindent
We must therefore use the two factors $(z-\ii)(z+\ii)=z^{2}+1$ in the
construction of the center series,

\noindent
\begin{eqnarray*}
  C_{z}
  & = &
  \frac{4}{\pi}\,
  \left(z^{2}+1\right)
  \sum_{j=0}^{\infty}
  \frac{(-1)^{j}}{2j+1}\,
  z^{2j+1}
  \\
  & = &
  \frac{4}{\pi}
  \left[
    \sum_{j=0}^{\infty}
    \frac{(-1)^{j}}{2j+1}\,
    z^{2j+3}
    +
    \sum_{j=0}^{\infty}
    \frac{(-1)^{j}}{2j+1}\,
    z^{2j+1}
  \right]
  \\
  & = &
  \frac{4}{\pi}
  \left[
    \sum_{j=1}^{\infty}
    \frac{(-1)^{j-1}}{2j-1}\,
    z^{2j+1}
    +
    z
    +
    \sum_{j=1}^{\infty}
    \frac{(-1)^{j}}{2j+1}\,
    z^{2j+1}
  \right]
  \\
  & = &
  \frac{4}{\pi}
  \left[
    z
    -
    \sum_{j=1}^{\infty}
    (-1)^{j}
    \left(
      \frac{1}{2j-1}
      -
      \frac{1}{2j+1}
    \right)
    z^{2j+1}
  \right]
  \\
  & = &
  \frac{4}{\pi}
  \left[
    z
    -
    \sum_{j=1}^{\infty}
    \frac{2(-1)^{j}}{4j^{2}-1}\,
    z^{2j+1}
  \right]
  \\
  & = &
  \frac{4}{\pi}\,
  z
  \left[
    1
    -
    \sum_{j=1}^{\infty}
    \frac{2(-1)^{j}}{4j^{2}-1}\,
    z^{2j}
  \right].
\end{eqnarray*}

\noindent
Unlike the original series, with coefficients that behave as $1/k$ (with
$k=2j+1$), this series has coefficients that go to zero as $1/k^{2}$ when
$k\to\infty$, and therefore is absolutely and uniformly convergent to a
continuous function. This shows, in particular, that our evaluation of the
set of dominant singularities of $w(z)$ was in fact correct. We have
therefore for $S_{z}$ the representation

\begin{displaymath}
  S_{z}
  =
  \frac{4}{\pi}\,
  \frac{z}{z^{2}+1}
  \left[
    1
    -
    \sum_{j=1}^{\infty}
    \frac{2(-1)^{j}}{4j^{2}-1}\,
    z^{2j}
  \right],
\end{displaymath}

\noindent
with the singularities factored out and where the series involved is
absolutely and uniformly convergent, and therefore converges much faster
than the original one.

We may now take the real and imaginary parts of the $S_{v}$ series in
order to obtain faster-converging representation of the original DP
Fourier series and its FC series. We have on the unit circle

\noindent
\begin{eqnarray*}
  \frac{z}{z^{2}+1}
  & = &
  \frac
  {z\left[(z^{*})^{2}+1\right]}
  {\left(z^{2}+1\right)\left[(z^{*})^{2}+1\right]}
  \\
  & = &
  \frac{z^{*}+z}{2+z^{2}+(z^{*})^{2}}
  \\
  & = &
  \frac{2\cos(\theta)}{2+2\cos(2\theta)}
  \\
  & = &
  \frac{\cos(\theta)}{1+\cos^{2}(\theta)-\sin^{2}(\theta)}
  \\
  & = &
  \frac{\cos(\theta)}{2\cos^{2}(\theta)}
  \\
  & = &
  \frac{1}{2\cos(\theta)},
\end{eqnarray*}

\noindent
and therefore

\noindent
\begin{eqnarray*}
  S_{v}
  & = &
  \frac{2}{\pi\cos(\theta)}
  \left[
    1
    -
    \sum_{j=1}^{\infty}
    \frac{2(-1)^{j}}{4j^{2}-1}\,
    \cos(2j\theta)
    -
    \ii
    \sum_{j=1}^{\infty}
    \frac{2(-1)^{j}}{4j^{2}-1}\,
    \sin(2j\theta)
  \right]
  \\
  & = &
  \frac{2}{\pi\cos(\theta)}
  \left[
    1
    -
    \sum_{j=1}^{\infty}
    \frac{2(-1)^{j}}{4j^{2}-1}\,
    \cos(2j\theta)
  \right]
  +
  \ii\,
  \frac{2}{\pi\cos(\theta)}
  \left[
    -
    \sum_{j=1}^{\infty}
    \frac{2(-1)^{j}}{4j^{2}-1}\,
    \sin(2j\theta)
  \right].
\end{eqnarray*}

\noindent
The original DP function is given by the real part,

\begin{displaymath}
  f_{\rm c}(\theta)
  =
  \frac{2}{\pi\cos(\theta)}
  \left[
    1
    -
    \sum_{j=1}^{\infty}
    \frac{2(-1)^{j}}{4j^{2}-1}\,
    \cos(2j\theta)
  \right],
\end{displaymath}

\noindent
and the corresponding FC function $f_{\rm s}(\theta)=\bar{f}_{\rm
  c}(\theta)$ is given by the imaginary part,

\begin{displaymath}
  f_{\rm s}(\theta)
  =
  -\,
  \frac{2}{\pi\cos(\theta)}
  \left[
    \sum_{j=1}^{\infty}
    \frac{2(-1)^{j}}{4j^{2}-1}\,
    \sin(2j\theta)
  \right].
\end{displaymath}

\noindent
Both of these series are absolutely and uniformly convergent.

\subsection{Another Regular Sine Series with Odd {\boldmath $k$}}

Consider the Fourier series of a periodic function built with segments of
quadratic functions, joined together so that the resulting function is
continuous and differentiable. As is well known it is given by the sine
series

\begin{displaymath}
  S_{\rm s}
  =
  \frac{32}{\pi^{3}}
  \sum_{j=0}^{\infty}
  \frac{1}{(2j+1)^{3}}\,
  \sin[(2j+1)\theta].
\end{displaymath}

\noindent
The corresponding FC series is then

\begin{displaymath}
  \bar{S}_{\rm s}
  =
  \frac{32}{\pi^{3}}
  \sum_{j=0}^{\infty}
  \frac{1}{(2j+1)^{3}}\,
  \cos[(2j+1)\theta].
\end{displaymath}

\noindent
Note that due to the factors of $1/k^{3}$ (with $k=2j+1$), these series
are already absolutely and uniformly convergent. But we will proceed with
the construction in any case. The complex $S_{v}$ series is given by

\begin{displaymath}
  S_{v}
  =
  \frac{32}{\pi^{3}}
  \sum_{j=0}^{\infty}
  \frac{1}{(2j+1)^{3}}\,
  v^{2j+1},
\end{displaymath}

\noindent
and the complex power series $S_{z}$ is given by

\begin{displaymath}
  S_{z}
  =
  \frac{32}{\pi^{3}}
  \sum_{j=0}^{\infty}
  \frac{1}{(2j+1)^{3}}\,
  z^{2j+1}.
\end{displaymath}

\noindent
The ratio test tells us that the disk of convergence of $S_{z}$ is the
unit disk. If we consider the inner analytic function $w(z)$ within this
disk we observe that $w(0)=0$, as expected. We have for this function

\begin{displaymath}
  w(z)
  =
  \frac{32}{\pi^{3}}
  \sum_{j=0}^{\infty}
  \frac{1}{(2j+1)^{3}}\,
  z^{2j+1}.
\end{displaymath}

\noindent
Being given by a monotonic series of step $2$ this function has two
dominant singularities, located at $z=1$ and at $z=-1$, as one can easily
verify by taking its second logarithmic derivative, which is proportional
to the inner analytic function of the standard square wave, that we
examined before in Subsection~\ref{APPregsqwave},

\noindent
\begin{eqnarray*}
  z\,
  \frac{dw(z)}{dz}
  & = &
  \frac{32}{\pi^{3}}
  \sum_{j=0}^{\infty}
  \frac{1}{(2j+1)^{2}}\,
  z^{2j+1},
  \\
  z\,
  \frac{d}{dz}
  \left[
    z\,
    \frac{dw(z)}{dz}
  \right]
  & = &
  \frac{32}{\pi^{3}}
  \sum_{j=0}^{\infty}
  \frac{1}{2j+1}\,
  z^{2j+1}.
\end{eqnarray*}

\noindent
We must therefore use the two factors $(z-1)(z+1)=z^{2}-1$ in the
construction of the center series,

\noindent
\begin{eqnarray*}
  C_{z}
  & = &
  \frac{32}{\pi^{3}}
  \left(z^{2}-1\right)
  \sum_{j=0}^{\infty}
  \frac{1}{(2j+1)^{3}}\,
  z^{2j+1}
  \\
  & = &
  \frac{32}{\pi^{3}}
  \left[
    \sum_{j=0}^{\infty}
    \frac{1}{(2j+1)^{3}}\,
    z^{2j+3}
    -
    \sum_{j=0}^{\infty}
    \frac{1}{(2j+1)^{3}}\,
    z^{2j+1}
  \right]
  \\
  & = &
  \frac{32}{\pi^{3}}
  \left[
    \sum_{j=1}^{\infty}
    \frac{1}{(2j-1)^{3}}\,
    z^{2j+1}
    -
    z
    -
    \sum_{j=1}^{\infty}
    \frac{1}{(2j+1)^{3}}\,
    z^{2j+1}
  \right]
  \\
  & = &
  \frac{32}{\pi^{3}}
  \left\{
    -z
    +
    \sum_{j=1}^{\infty}
    \left[
      \frac{1}{(2j-1)^{3}}
      -
      \frac{1}{(2j+1)^{3}}
    \right]
    z^{2j+1}
  \right\}
  \\
  & = &
  \frac{32}{\pi^{3}}
  \left[
    -z
    +
    \sum_{j=1}^{\infty}
    \frac{24j^{2}+2}{\left(4j^{2}-1\right)^{3}}\,
    z^{2j+1}
  \right]
  \\
  & = &
  \frac{32}{\pi^{3}}\,
  z
  \left[
    -1
    +
    \sum_{j=1}^{\infty}
    \frac{24j^{2}+2}{\left(4j^{2}-1\right)^{3}}\,
    z^{2j}
  \right].
\end{eqnarray*}

\noindent
Unlike the original series, with coefficients that behave as $1/k^{3}$
(with $k=2j+1$), this series has coefficients that go to zero as $1/k^{4}$
when $k\to\infty$, and therefore converges faster than the original one.
This shows, in particular, that our evaluation of the set of dominant
singularities of $w(z)$ was in fact correct. We have therefore for $S_{z}$
the representation

\begin{displaymath}
  S_{z}
  =
  \frac{32}{\pi^{3}}\,
  \frac{z}{z^{2}-1}
  \left[
    -1
    +
    \sum_{j=1}^{\infty}
    \frac{24j^{2}+2}{\left(4j^{2}-1\right)^{3}}\,
    z^{2j}
  \right],
\end{displaymath}

\noindent
with the singularities factored out. Although both this series and the
original one are absolutely and uniformly convergent, this converges
faster, and may be differentiated twice, still resulting in other series
that are also absolutely and uniformly convergent.

We may now take the real and imaginary parts of the $S_{v}$ series in
order to obtain faster-converging representation of the original DP
Fourier series and its FC series. We have on the unit circle, as we saw
before in Subsection~\ref{APPregsqwave},

\begin{displaymath}
  \frac{z}{z^{2}-1}
  =
  \frac{-\ii}{2\sin(\theta)},
\end{displaymath}

\noindent
and therefore

\noindent
\begin{eqnarray*}
  S_{v}
  & = &
  \frac{16}{\pi^{3}}\,
  \frac{-\ii}{\sin(\theta)}
  \left[
    -1
    +
    \sum_{j=1}^{\infty}
    \frac{24j^{2}+2}{\left(4j^{2}-1\right)^{3}}\,
    \cos(2j\theta)
    +
    \ii
    \sum_{j=1}^{\infty}
    \frac{24j^{2}+2}{\left(4j^{2}-1\right)^{3}}\,
    \sin(2j\theta)
  \right]
  \\
  & = &
  \frac{16}{\pi^{3}\sin(\theta)}
  \left[
    \sum_{j=1}^{\infty}
    \frac{24j^{2}+2}{\left(4j^{2}-1\right)^{3}}\,
    \sin(2j\theta)
  \right]
  +
  \\
  &   &
  +
  \ii\,
  \frac{16}{\pi^{3}\sin(\theta)}
  \left[
    1
    -
    \sum_{j=1}^{\infty}
    \frac{24j^{2}+2}{\left(4j^{2}-1\right)^{3}}\,
    \cos(2j\theta)
  \right].
\end{eqnarray*}

\noindent
The original DP function is given by the imaginary part,

\begin{displaymath}
  f_{\rm s}(\theta)
  =
  \frac{16}{\pi^{3}\sin(\theta)}
  \left[
    1
    -
    \sum_{j=1}^{\infty}
    \frac{24j^{2}+2}{\left(4j^{2}-1\right)^{3}}\,
    \cos(2j\theta)
  \right],
\end{displaymath}

\noindent
and the corresponding FC function $f_{\rm c}(\theta)=\bar{f}_{\rm
  s}(\theta)$ is given by the real part,

\begin{displaymath}
  f_{\rm c}(\theta)
  =
  \frac{16}{\pi^{3}\sin(\theta)}
  \left[
    \sum_{j=1}^{\infty}
    \frac{24j^{2}+2}{\left(4j^{2}-1\right)^{3}}\,
    \sin(2j\theta)
  \right].
\end{displaymath}

\noindent
Both of these series converge faster than the original Fourier series.

\end{document}